\documentclass[a4paper,12pt,twoside]{amsart}
\usepackage{amsmath,amsthm,amsfonts}
\usepackage{amssymb}
\usepackage[T1]{fontenc}
\usepackage{graphics,graphicx}

\newtheorem{thm}{Theorem}[section]

\newtheorem{cor}[thm]{Corollary}
\newtheorem{lem}[thm]{Lemma}

\newtheorem{prop}[thm]{Proposition}
\newtheorem{proposition}[thm]{Proposition}

\newtheorem{defns}[thm]{Definitions}

\theoremstyle{definition}

\newtheorem{rem}[thm]{Remark}

\newcommand{\cal} {\mathcal}

\newcommand{\R}{{\mathbb{R}}}

\newcommand{\E}{{\mathbb{E}}}
\newcommand{\Q}{{\mathbb{Q}}}
\newcommand{\C}{{\mathbb{C}}}

\newcommand{\Z}{{\mathbb{Z}}}
\newcommand{\N}{{\mathbb{N}}}
\newcommand{\T}{{\mathbb{T}}}
\newcommand{\PP}{{\mathbb{P}}}

\def \RR {\R}

\def \EE {\E}
\def \QQ {\Q}
\def \CC {\C}

\def \ZZ {\Z}
\def \NN {\N}
\def \TT {\T}

\newcommand{\SL}{\operatorname{SL}}
\newcommand{\GL}{\operatorname{GL}}
\newcommand{\Ree}{\operatorname{\Re e}}

\newcommand{\supp}{\operatorname{supp}}

\newcommand{\Id}{\operatorname{Id}}
\newcommand{\Aut}{\operatorname{Aut}}
\newcommand{\AP}{\operatorname{AP}}

\def\P{{ \bold P}}

\def\eop{\qed}
\def\Proof {\vskip -2mm {{\it Proof}.}}
\def\proof {\vskip -2mm {{\it Proof}.}}

\def\trace{{\rm trace}}
\def\dim{{\rm dim}}
\def\det{{\rm det}}

\font\fivegoth=eufm5 \font\sevengoth=eufm7 \font\tengoth=eufm10
\newfam\gothfam
\scriptscriptfont\gothfam=\fivegoth \textfont\gothfam=\tengoth
\scriptfont\gothfam=\sevengoth

\title [Ergodicity of group actions, spectral gap and applications]
{Ergodicity of group actions and spectral gap, applications to
random walks and Markov shifts}

\medskip
\author{J.-P. Conze}
\author{Y. Guivarc'h}
\address{J.-P. Conze, Y. Guivarc'h,
\hfill \break
IRMAR, CNRS UMR 6625
\hfill \break
Universit\'e de Rennes 1
\hfill \break
Campus de Beaulieu, 35042 Rennes Cedex, France}
\email{
\hfill \break
conze@univ-rennes1.fr
\hfill \break
guivarch@univ-rennes1.fr}

\date{May 2011}


\begin{document}

\begin{abstract} Let $(X, \cal B, \nu)$ be a probability space and let $\Gamma$
be a countable group of $\nu$-preserving invertible maps of $X$ into
itself. To a probability measure $\mu$ on $\Gamma$ corresponds a
random walk on $X$ with Markov operator $P$ given by $P\psi(x) =
\sum_{a} \psi(ax) \, \mu(a)$. A powerful tool is the spectral gap
property for the operator $P$ when it holds. We consider various
examples of ergodic $\Gamma$-actions and random walks and their
extensions by a vector space: groups of automorphisms or affine
transformations on compact nilmanifolds, random walk in random
scenery on non amenable groups, translations on homogeneous spaces
of simple Lie groups, random walks on motion groups. The spectral
gap property is applied to obtain limit theorems,
recurrence/transience property and ergodicity for random walks on
non compact extensions of the corresponding dynamical systems.
\end{abstract}

\maketitle

\tableofcontents

\section*{Introduction}

Let $(X, \cal B, \nu)$ be a metric space endowed with its Borel $\sigma$-algebra
$\cal B$ and a probability measure $\nu$, and let $\Gamma$ be a
countable group of Borel invertible maps of $X$ into itself which
preserve $\nu$.

Let $\mu$ be a probability measure on $\Gamma$ such that the group
generated by $A:=\supp(\mu)$ is $\Gamma$. We consider the random
walk on $X$ defined by $\mu$, with Markov operator $P$ given by
\begin{eqnarray}
P\psi(x) &=& \sum_{a \in A} \psi(ax) \, \mu(a), \ x \in X.
\label{MarkovP}
\end{eqnarray}

These data, i.e., the probability space $(X, \nu)$, the group
$\Gamma$ acting on $(X, \nu)$ and the probability measure $\mu$ on
$\Gamma$, will be denoted by $(X, \nu, \Gamma, \mu)$.

The operator $P$ is a contraction of $L^p(X, \nu)$, $ \forall p\geq
1$, and it preserves the subspace $L_0^2(X, \nu)$ of functions
$\varphi$ in $L^2(X, \nu)$ such that $\nu(\varphi) = 0$. $P$ is said
to be {\it ergodic} if the constant functions are the only
$P$-invariant functions in $L^2(X, \nu)$.

Ergodicity of $P$ is equivalent to ergodicity of the action
of $\Gamma$ on the measure space $(X, \cal B,\nu)$. Indeed, any
$\Gamma$-invariant function is obviously $P$-invariant. Conversely,
if $\varphi$ in $L^2(\nu)$ is $P$-invariant, then, by strict
convexity of $L^2(X, \nu)$, we have $\varphi(ax) = \varphi(x)$,
$\nu$-a.e. for every $a \in \supp(\mu)$. Therefore $\varphi$ is
$\Gamma$-invariant, hence $\nu$-a.e. constant if $\Gamma$ acts
ergodically on $(X, \nu)$.

Our aim is to consider some examples of ergodic actions and
extensions of these actions by a vector space. We will use a strong
reinforcement of the ergodicity, the spectral gap property for the
operator $P$ when it holds and we will develop some of its
consequences. Let us recall its definition and related notions.

\begin{defns} \label{def01} {\rm We denote by $\Pi_0$ the restriction of $P$ defined by
(\ref{MarkovP}) to $L_0^2(X, \nu)$ and by $r(\Pi_0):=\lim_n
\|\Pi_0^n\|^{1\over n}$ its spectral radius. If $r(\Pi_0) < 1$, we
say that $(X, \nu, \Gamma, \mu)$ satisfies the {\it spectral gap
property} (we will use the shorthand "property (SG)").}
\end{defns}

We recall that a unitary representation $\rho$ of a group $\Gamma$
in a Hilbert space $H$ is said to {\it contain weakly} the identity
representation if there exists a sequence $(x_n)$ in $H$ with
$\|x_n\| = 1$ such that, for every $\gamma \in \Gamma$, $\lim_n
\|\rho(\gamma) x_n - x_n \|= 0$. See \cite{BeHaVa08} for this
notion.

Recall also that $\Gamma$ is said to have property (T) if, when the
identity representation is weakly contained in a unitary
representation $\rho$ of $\Gamma$, then it is contained in $\rho$.

The natural action of $\Gamma$ on $L_0^2(X)$ defines a unitary
representation $\rho_0$ of $\Gamma$ in $L_0^2(X)$. Property (SG)
implies that the identity representation of $\Gamma$ is not weakly
contained in $\rho_0$. The converse is true if, for every $k
> 0$, $(\supp (\mu))^k$ generates $\Gamma$ (see below Corollary
\ref{cor312}). Property (SG) depends only on the support of $\mu$.

For a countable group $\Gamma$ acting measurably on a probability
measure space $(X,\nu)$ where $\nu$ is $\Gamma$-invariant, according
to \cite{FuSh99} the $\Gamma$-action on $(X,\nu)$ is said to be {\it
strongly ergodic} if $\nu$ is the unique $\Gamma$-invariant
continuous positive normalized functional on $L^\infty(X,\nu)$.
Property (SG) implies strong ergodicity, hence ergodicity of the
action of $\Gamma$ on $X$.

Our framework will be essentially algebraic. As examples, we study
the action of groups of automorphisms or affine transformations on
tori and compact nilmanifolds, and translations on homogeneous
spaces of simple Lie groups. In Section \ref{nilpSection} we show
for nilmanifolds that the ergodicity of $P$ follows from the
ergodicity of its restriction to the maximal torus quotient. In
Section \ref{SGProp}, we recall property (SG) for subgroups of
$\SL(d, \ZZ)$ acting on $\TT^d$, as well as recent results on
property (SG) for the nilmanifolds. In Section \ref{appli} we
consider random walks on non compact extensions of dynamical systems
and apply property (SG) to recurrence and ergodicity. The last
section is devoted to examples.

The authors thank Bachir Bekka for useful discussions.

 \vskip 5mm
\section{Ergodicity of a group of affine transformations on nilmanifolds}
\label{nilpSection}

In this section, we consider groups of affine transformations
$\Gamma$ on compact nilmanifolds $X$. In order to obtain ergodicity
of Markov operators on $X$, as described in the introduction, we
study the question of ergodicity of the action of $\Gamma$.

Let $N$ be a connected, simply connected, nilpotent Lie group and
$D$ a lattice in $N$, i.e. a discrete subgroup $D$ such that the
quotient $X = N/D$ is compact. If $L_1, L_2$ are two subgroups of
$N$, we denote by $[L_1, L_2]$ the closed subgroup generated by the
elements $\{n_1\,n_2\,n_1^{-1} \, n_2^{-1}, n_1 \in L_1, n_2 \in L_2
\}$, $L':=[L,L]$ the derived group of $L$, $e$ the neutral element of $N$.
The descending series of $N$ is defined by
$$N \supset N^1\supset ...\supset N^{k-1} \supset N^{k}
\supset \{e\},$$ where $N^{\ell+1}:= [N^\ell, N]$, for $\ell \geq
0$, with $N^0 = N$.

The elements $g \in N$ act on $N/D$ by left translation: $nD \in N/D
\rightarrow gnD$. We say that $\tau$ is an automorphism of the
nilmanifold $N/D$ if $\tau$ is an automorphism of the group $N$ such
that $\tau D = D$. The group of automorphisms of $N/D$ is denoted by
$\Aut(N/D)$. The action of $\tau \in \Aut(N/D)$ on $N/D$ is $nD
\rightarrow \tau(n)D$. An affine transformation $\gamma$ of $N/D$ is
a map of the form:
\begin{equation} nD \rightarrow \gamma(n)D = \alpha_\gamma \,
\tau_\gamma(n) D,\label{trans-aff}
\end{equation}
with $\alpha_\gamma \in N$ and $\tau_\gamma \in \Aut(N/D)$.

Let $\Gamma$ be a group of affine transformations of the
nilmanifold. The measure $m$ on $N/D$ deduced from a Haar measure on
$N$ is $\Gamma$ invariant. The group $\Gamma$ acts on the quotients
$N^\ell/N^\ell \cap D$, $0\leq \ell \leq k+1$, and in particular on
the factor torus $T= N/N^1.D$. When $\Gamma$ is a group of automorphisms,
ergodicity of the action on the torus is equivalent
to the fact that every non trivial character has an infinite $\Gamma$-orbit.

When $\Gamma$ is generated by a single
automorphism (or more generally by an affine transformation), W.
Parry has proved (\cite{Pa69}, \cite{Pa70}) that the ergodicity of
the action on the quotient $T$ implies the ergodicity of the action on
the nilmanifold. We will show (Theorem \ref{main}) that the analogous statement holds
for a group of affine transformations.

{\bf Notations}: For a given group $\Gamma$ of affine
transformations of $N/D$, $\tilde \Gamma$ denotes the subgroup
of $\Aut(N/D)$ generated by $\{\tau_\gamma, \gamma \in \Gamma \}$,
where $\tau_\gamma$ is the automorphism associated with $\gamma$ as in
(\ref{trans-aff}). We denote by $N_e^\ell$ the Lie algebra of $N^\ell$ and by $d\tau_e$
the linear map tangent at $e$ to an automorphism $\tau$ of $N$.

We will use the following lemmas.
\begin{lem} \label{distal} (cf. {CoGu74}) If $\Gamma$ is a subgroup
of $\GL(\RR^d)$ such that the eigenvalues of each element of $\Gamma$ has
modulus 1, then there is a $\Gamma$-invariant subspace $W \not = \{0\}$
of $\RR^d$ such that the action of $\Gamma$ on $W$ is relatively compact.
If $\Gamma$ is a subgroup of $\GL(\ZZ^d)$, the action of $\Gamma$ on $W$
is that of a finite group of rotations and reduces to the identity for
$\gamma$ in a subgroup $\Gamma_0$ of finite index in $G$.
\end{lem} \Proof \ \ We extend the action of $\Gamma$ to $\CC^d$.
Let $\tilde W$ be a subspace of $\CC^d$  which is different from $\{0\}$
and invariant by $\Gamma$ on which the action of $\Gamma$ is
irreducible. Let  $(e_i)$ be a basis of $\tilde W$, and let $E_{ij}$
be the maps defined by $E_{ij}(e_k) = \delta_{kj} e_i, \forall k$.
We denote by $\tilde \tau$ the endomorphism corresponding to the
action of $\tau \in \Gamma$ on $\tilde W$.

The trace of each automorphism $\tilde \tau$, for $\tau \in \Gamma$,
satisfies: $\trace(\tilde \tau) \leq C \dim(\tilde W)$.

The action of $\Gamma$ on $\tilde W$  being irreducible, by
Burnside's theorem there are constants $b_k$
and elements $\tau_k$ of $\Gamma$ such that $E_{ji} = \sum_k b_k
\tilde \tau_k$. The coefficients of the transformations $\tilde \tau$
satisfy then:
$$|a_{ij}(\tau)| = |\trace(\tilde \tau E_{ji})|\leq \sum_k |b_k|
|\trace(\tilde \tau \tilde \tau_k)| \leq C \dim(\tilde W) \sum_k
|b_k|.$$ Therefore $\sup_{\tau \in \Gamma} |a_{i j}(\tau| < \infty$,
which implies the relative compactness of the action of $\Gamma$ on $\tilde W$,
as well on $W$, the $\Gamma$-invariant subspace of $\RR^d$ generated
by $\{ \Ree v, v \in \tilde W \}$.

Now assume that $\Gamma$ is a subgroup of $\GL(\ZZ^d)$. The symmetric
functions of the eigenvalues of $\gamma$ in $\Gamma$ take values in $\ZZ$ and
remain bounded when $\gamma$ runs in $\Gamma$. This
implies that the set of the characteristic polynomials of the
elements $\gamma$ is finite. If $\lambda$ is an eigenvalue of
$\gamma$ the set $(\lambda^n)_{n \in \ZZ}$ is finite and therefore
$\lambda$ is a root of the unity. The order of these roots remains
bounded on $\Gamma$. This implies the last assertion. \eop

\begin{lem} \label{subgroup1} If a group $\Gamma$ of affine transformations
on $\TT^d$ has an invariant square integrable function $f$ non a.e.
constant, then it has an invariant function which is a non
identically constant trigonometric polynomial. If the action of
$\Gamma$ is ergodic, every eigenfunction is a trigonometric
polynomial.
\end{lem} \Proof \ \ Let $f \in L^2(N/D)$ be a $\Gamma$-eigenfunction, $f\circ
\gamma = \beta(\gamma) f, \forall \gamma \in \Gamma$. We have, for
every $\gamma \in \Gamma$~:
\begin{equation} f = \sum_{p \in \ZZ^d} \hat f(p) e^{2\pi i < p,.>} =
\overline{\beta(\gamma)} \sum_{p \in \ZZ^d} \hat f(p) e^{2\pi i <
p,\alpha_\gamma>} e^{2\pi i <{}^t\tau_\gamma p,.>};
\end{equation}
hence: $|\hat f(p)| = |\hat f({}^t\tau p)|, \forall p\in \ZZ^d$.

Let $R:=\{p \in \ZZ^d: |\hat f(p)| \not = 0 \}$. For two
automorphism $\tau, \tau'$ of the torus and $p \in \ZZ^d$ such that
${}^t\tau p \not = {}^t\tau' p$, the characters $e^{2\pi i <{}^t\tau
p,.>}$ and $e^{2\pi i <{}^t\tau p,.>}$ are orthogonal. Therefore the
orbit $\{{}^t\tau_\gamma p,\  \gamma \in \Gamma\}$ of every element
$p$ of $R$ is finite. The set $R$ decomposes into finite disjoint
subsets $R_k$, with each $R_k$ permuted by the automorphisms
$\tau_\gamma \in \tilde \Gamma$.

The subspaces $W_k$ of $L^2$ generated by $e^{2 \pi i <p, .>}$, for
$p \in R_k$, have a finite dimension, are pairwise orthogonal and
invariant by each $\gamma \in \Gamma$. The orthogonal
projections of $f$ on these subspaces give $\Gamma$-eigenfunctions
with the same eigenvalue as for $f$. This shows the existence
of a non constant eigenfunction (invariant if $f$ is invariant)
which is a trigonometric polynomial.
If the group $\Gamma$ acts ergodically, only one of these
projections is non null. Hence $f$ is a trigonometric polynomial.
\eop

\begin{lem} \label{subgroup2} If a group of affine transformations
$\Gamma$ of a torus $\TT^d$ is ergodic, then
every subgroup $\Gamma_0$ of $\Gamma$ with finite index is also
ergodic on $\TT^d$.
\end{lem} \Proof \ \ Let $\Gamma_0$ be a subgroup of $\Gamma$
with finite index. As the action of $\Gamma$ is ergodic, the
$\sigma$-algebra of the $\Gamma_0$-invariant subsets is an atomic
finite $\sigma$-algebra  whose elements are permuted by $\gamma \in
\Gamma$. From Lemma \ref{subgroup1} there exists a non constant
trigonometric polynomial which is invariant by $\Gamma_0$. This
polynomial should be measurable with respect to the $\sigma$-algebra
of the $\Gamma_0$-invariant subsets which is atomic. The
connectedness of the torus implies that it is constant. \eop

\goodbreak
\vskip 3mm {\bf Ergodicity of a group of affine transformations}
\begin{thm} \label{main} Let $\Gamma$ be a group of affine transformations
on $N/D$. If its action on the torus quotient $N/N^1.D$ is ergodic,
then every eigenfunction for the action of $\Gamma$ on $N/D$
factorizes into an eigenfunction on $N/N^1.D$. In particular, the
action of $\Gamma$ is ergodic on $N/D$ if and only if its action on
the quotient $N/N^1.D$ is ergodic.
\end{thm}
\Proof  \ \ We follow essentially the method of W. Parry
(\cite{Pa69}). We make an induction on the length $k$ of the
descending central series of $N$. The property stated in the theorem
is clearly satisfied if $k= 0$.

The induction assumption is that, for every group of affine
transformations of $N/D$, ergodicity of the action on $N/N^1.D$
implies ergodicity of the action on $N/N^k.D$ and every
eigenfunction for the action on $N/N^k.D$ factorizes through
$N/N^1.D$. (The quotient $(N/N^k)/(N/N^k)'$ can be identified with the
quotient $N/N^1$.)

Remark that, for every subgroup $\Gamma_0$ with a finite index in a
group $\Gamma$ of affine transformations on $N/D$, the action of
$\Gamma_0$ on $N/N^1.D$ is ergodic if the action of $\Gamma$ on
$N/N^1.D$ is also ergodic from Lemma \ref{subgroup2}. The induction
assumption implies then that $\Gamma_0$ acts ergodically on
$N/N^k.D$ and that every eigenfunction for the action of $\Gamma_0$
on $N/N^k.D$ factorizes through $N/N^1.D$.

Let $f \in L^2(N/D)$ be a $\Gamma$-eigenfunction, i.e. such that for
complex numbers $\beta(\gamma)$ of modulus 1,
\begin{equation} f(\gamma(n)D) = f(\alpha_\gamma\, \tau_\gamma(n) D)
= \beta(\gamma) f(nD), \forall \gamma \in \Gamma. \label{fct-propre}
\end{equation}

We are going to show that $f$ factorizes into an eigenfunction on
the quotient $N/ N^k.D$ and therefore, by the induction hypothesis,
into an eigenfunction on $N/N^1.D$.

The proof is given in several steps.

a) We denote by $Z$ the center of $N$. We have $N^k \subset Z\cap
N^{k-1}$. The torus $H:= Z \cap N^{k-1}/ Z \cap N^{k-1} \cap D$ acts
by left translation on $N/D$ and its action commutes with the
translation by elements of $N$. Let $\Theta$ be the group of
characters of $H$. The space $L^2(N/D)$ decomposes into pairwise
orthogonal subspaces $V_\theta$, where $\theta$ belongs to $\Theta$
and $V_\theta$ stands for the subspace of functions which are
transformed according to the character $\theta$ under the action of
$H$.
$$V_\theta = \{\varphi \in L^2(N/D) : \varphi(h nD)
= \theta(h)\  \varphi(nD), \forall h \in H\}.$$

If $\tau$ is an automorphism of $N/D$, $h \in H
\rightarrow \theta(\tau (h))$ defines a character on $H$ denoted by
$\tau \theta$. We have $\varphi \in V_\theta \Leftrightarrow \varphi
\circ \gamma \in V_{\tau_\gamma \theta}$.

We decompose $f$ into components in the subspaces $V_\theta$. By
(\ref{fct-propre}), we have, for every $\gamma \in \Gamma$, the
orthogonal decomposition of $f$:
\begin{equation} f = \sum_{\theta \in \Theta} f_\theta =
\overline{\beta(\gamma)} \sum_{\theta \in \Theta} \ \ f_\theta\circ
\gamma,\label{decomp}
\end{equation}
with $f_\theta \in V_{\theta}$, $f_\theta\circ \gamma \in
V_{\tau_\gamma \theta}$.  We will show that the components
$f_\theta$, hence also $f$, are invariant by translation by the
elements of $N^k$.

Let us fix $\theta$ such that $\|f_\theta \|_2 \not = 0$. Let $\Gamma_0
:= \{\gamma \in \Gamma : \tau_\gamma \theta = \theta \}$. For two
automorphisms $\tau, \tau'$ of $N/D$, if $\tau \theta \not = \tau'
\theta$, the subspaces $V_{\tau \theta}$ and $V_{\tau' \theta}$ are
orthogonal. The equality of the norms $\|f_\theta \circ\gamma \|_2 =
\|f_\theta \|_2, \forall \gamma \in \Gamma$ and Equation
(\ref{decomp}) imply that there are only a finite number of distinct
images in the orbit $\{\tau_\gamma \theta, \gamma \in \Gamma\}$. The
group $\Gamma_0$ has a finite index in $\Gamma$. As remarked above,
its action on $N/ N^k.D$, as the action of $\Gamma$, is ergodic.

\vskip 3mm From now on, we consider $\Gamma_0$ and the component
$f_\theta$ which we denote by $f$ for simplicity.

From what precedes, we have an ergodic action of a group of affine
transformations $\Gamma_0$, a character $\theta \in \Theta$ such
that $\tau_\gamma \theta = \theta, \forall \gamma \in \Gamma_0$, and
a function $f$ satisfying (denoting by $x$ an element of $N/D$):
$$f(\alpha_\gamma \tau_\gamma(x)) = \beta(\gamma) f(x), \ f( h.x) =
\theta(h) f(x), \forall h \in Z\cap N^{k-1}.$$

We can assume that $\theta$ is non trivial on $N^k$. By replacing
$N/D$ by $N/H_0. D$, where $H_0$ is the connected component of the
neutral element of $Ker \, \theta$, we can also assume that $N^k /
N^k \cap D$ has dimension 1.

The previous equations imply that $|f|$ is $\gamma$-invariant for
every $\gamma \in \Gamma_0$ and $N^k$-invariant. Therefore the
function $|f|$ is a.e. equal to a constant that we can assume to be
1.

\vskip 3mm b) For $g \in N$, let
$\theta(g) := \int f(gx) \overline{f(x)} \, dx$.
The function $\theta$ is continuous and $\theta(e) = 1$;
therefore $\theta(g) \not = 0$ on a neighborhood of $e$. The invariance of the measure implies:
\begin{equation}
\theta(g^{-1}) = \overline{\theta(g)}. \label{conjug}
\end{equation}

Denote by $G$ the subgroup of $N^{k-1}$ defined by
\begin{equation} G := \{g \in N^{k-1}: f(g x) = \theta(g) f(x) \}
= \{g \in N^{k-1}: |\theta(g)| = 1 \},\label{defG}
\end{equation}
(the equality in (\ref{defG}) follows from the equality in
Cauchy-Schwarz inequality). Denote by $G_0$ the connected component
of the neutral element in $G$.

For $g$ in $N^{k-1}$ and $h$ in $N$, we have: $\theta(hgh^{-1}) =
\theta(hgh^{-1}g^{-1}) \theta(g)$; hence
\begin{equation}
|\theta(hgh^{-1})| = |\theta(g)|, \ \forall  g \in N^{k-1}, h \in N.
\label{modulgh}
\end{equation}

For $g$ in $G$ and $h$ in $N$, the relation $f(g x) = \theta(g)
f(x)$ implies $\theta(gh) = \theta(g)\, \theta(h)$. Therefore we
have by applying (\ref{conjug}) to $gh$~:
\begin{equation}
g \in G, h \in N \Rightarrow \theta(gh) = \theta(hg) = \theta(g)\,
\theta(h). \label{thetagh}
\end{equation}

Equation (\ref{modulgh}) implies: $|\theta(hgh^{-1})| = 1, \ \forall
g \in G, h \in N$. The group $G$ (and therefore $G_0$) is a closed
normal subgroup of $N$.

From the equation of eigenfunction (\ref{fct-propre}), we
have:
\begin{eqnarray} f(\alpha_\gamma \tau_\gamma (g) \, \tau_\gamma (x))
\overline{f(\alpha_\gamma \tau_\gamma (x))} &=& f(\gamma (g x))
\overline{f(\gamma (x))} \nonumber\\
&=& \beta(\gamma) \overline{\beta(\gamma)} f(gx) \overline{f(x)} =
f(gx) \overline{f(x)} \label{prodf}\end{eqnarray}and therefore, by
invariance of the measure:
\begin{eqnarray*}
\int f(\alpha_\gamma \tau_\gamma(g) \alpha_\gamma^{-1} x)
\overline{f(x)} \ dx = \int f(\alpha_\gamma \tau_\gamma(g)
\tau_\gamma(x)) \overline{f(\alpha_\gamma \tau_\gamma(x))} \ dx =
\int f(g x) \overline{f(x)} \ dx,
\end{eqnarray*}
which implies: $\theta(\alpha_\gamma \tau_\gamma(g)
\alpha_\gamma^{-1}) = \theta(g)$; hence, from (\ref{modulgh})~:
\begin{equation}|\theta(\tau_\gamma g)| = |\theta(g)|, \forall
\gamma \in \Gamma_0.\label{moduTheta} \end{equation}

We define two subsets $\exp W_1$ and $\exp W_2$ containing
respectively the stable and instable subgroups of the automorphisms
$\tau$ in $\tilde \Gamma_0$ acting on $N^{k-1}/N^k$ by setting
\begin{eqnarray*} W_1 &=& \{v \in N_e^{k-1} : \exists \tau \in
\tilde \Gamma_0 : \lim_{n
\rightarrow +\infty} \ d\tau_e^n v \ {\rm mod \ } N_e^k = 0 \}, \\
W_2 &=& \{v \in N_e^{k-1} :  \exists \tau \in \tilde \Gamma_0 :
\lim_{n \rightarrow -\infty} \ d\tau_e^n v \ {\rm mod \ } N_e^k = 0
\}. \end{eqnarray*} Let us show that $\exp \,W_i \subset G_0$,
$i=1,2$. Let $g \in \exp W_1$. It belongs to a one parameter subgroup $(g_t)$
such that, for every $t$, there exists a sequence  $(g_n)$ of
elements of $N^k$ such that $\lim_n \tau^n (g_t) \, g_n = e$. From
(\ref{moduTheta}), this implies: $$|\theta(g_t)| = |\theta(\tau^n
(g_t)  \, g_n)| = |\theta(\tau^n (g_t)|  \, |\theta(g_n)|  =
|\theta(\tau^n (g_t)| \rightarrow  |\theta(e)| = 1.$$ Therefore
$g_t$ is in $G$, for every $t$, hence $g \in G_0$.
The analysis is the same for $W_2$.

\vskip 3mm c) Now we prove: $[N, G_0] = N^k$.

As we are reduced to the case where $N^k$ is of dimension
1, the other possibility is that $[N, G_0] = \{e\}$. Let us assume
that $[N, G_0] = \{e\}$. The subgroup $G_0$ is then in the center of
$N$ and contained in $Z \cap N^{k-1}$.

Consider the quotient $N^{k-1} / Z \cap N^{k-1}$. As $G_0$ contains
$\exp \,W_1$ and $\exp \,W_2$, the differentials of the
automorphisms $\tau \in \Gamma_0$ have only eigenvalues of modulus 1
for their action on $N_e^{k-1} / Z_e \cap N_e^{k-1}$.

By Lemma \ref{distal}, there exists then in $N_e^{k-1} / Z_e \cap
N_e^{k-1}$ a subspace $W_3$ on which the action of the
transformations which are linear tangent to the automorphisms $\in
\Gamma_0$ is compact. As the automorphisms $\tau$ preserve a
lattice, their action on $W_3$ has a finite order. The subgroup
$\tilde \Gamma_1$ which leaves fixed the elements of $W_3$ has a
finite index in $\tilde \Gamma_0$.

Let $g$ be an element of $N^{k-1}$ in $\exp \, W_3$. For every
$\tau_\gamma \in \tilde \Gamma_1$, there exists $g_0 \in Z\cap
N^{k-1}$ such that $\tau_\gamma g = g_0g$. Equation (\ref{prodf}) reads
\begin{eqnarray*}
f(\alpha_\gamma \, g_0 g \, \tau_\gamma (x))
\overline{f(\alpha_\gamma \tau_\gamma (x))} &=& f(g_0 \alpha_\gamma
g \alpha_\gamma^{-1} g^{-1} \, g \, \alpha_\gamma \,\tau_\gamma (x))
\overline{f(\alpha_\gamma \tau_\gamma (x))} \\
&=& \theta(g_0) \, \theta(\alpha_\gamma  g \alpha_\gamma^{-1}
g^{-1}) f(g\,\gamma(x)) \overline{f(\gamma(x))} = f(gx)
\overline{f(x)}.
\end{eqnarray*}
The function $f(g.) \ \overline{f(.)}$ is an eigenfunction for every
$\gamma  \in \Gamma_1$ and is invariant by $h \in N^k$. It
factorizes into an eigenfunction for $\Gamma_1$ on $N/N^k.D$. Either
its integral is 0 (if for $\gamma \in \Gamma_1$ the corresponding
eigenvalue is $\not = 1$) or it is invariant by $\Gamma_1$. In the
latter case, as by the induction hypothesis ergodicity holds for the
action of $\Gamma_1$ on $N/N^k.D$, the function $f(g.)  \
\overline{f(.)}$ is equal to a constant which has modulus 1 (since
$|f| = 1$).

We have therefore $|\theta(g)| = 0$ or 1. By a continuity argument
$|\theta(g)| = 1$, which implies that $g \in G$. Using as above a
one parameter subgroup, we obtain that $g \in G_0$. This gives a
contradiction, since $G_0 \subset Z \cap N^{k-1}$.

\vskip 3mm d) Let $(h_t)$ be a one parameter subgroup of $N$ and let $g$ be in
$G_0$. We have from (\ref{thetagh}):
\begin{eqnarray*}
\theta(h_t) \, \theta(g) &=& \theta(h_t g) = \theta(h_t g h_t^{-1} g^{-1} g h_t) \\
&=&\theta(h_t g h_t^{-1} g^{-1}) \theta(g h_t) =\theta(h_t g h_t^{-1}
g^{-1}) \theta(g) \, \theta(h_t).
\end{eqnarray*}

By continuity, $\theta(h_t)$ is different from zero in a
neighborhood of $t=0$. The previous relation implies
$\theta(h_t g h_t^{-1} g^{-1}) = 1$ in a neighborhood of
$t=0$ and, $G_0$ being connected, is equal to 1 everywhere. As $[N,
G_0] = N^k$, this shows that the character $\theta$ is identically
equal to 1 on $N^k$ and therefore the announced factorization property
is satisfied. \eop

{\it Remark} There are compact nilmanifolds $N/\Gamma$ for which
the group $\Aut(N/\Gamma)$ is non ergodic (cf. \cite{DiLi57}).
This contrasts with the case of Heisenberg nilmanifolds, for which there is
a large group of automorphisms.

\vskip 3mm {\bf An example}

Now we give an example of nilmanifold with an ergodic group $\Gamma$
of automorphisms such that each automorphism in $\Gamma$ is non ergodic.

\vskip 3mm{\it Construction on the torus}

Examples of groups of matrices such that each of them has an eigenvalue
equal to 1 can be constructed by action on the space of quadratic forms.
We explicit the example in dimension 2.

Let $A := \begin{pmatrix} a & b \\ c & d \end{pmatrix} \in
\GL(2,\RR)$, with eigenvalues $\lambda_1, \lambda_2$. The matrix corresponding
to the action of $A$ on the vector space of symmetric $2\times2$ matrices
$M \rightarrow A M A^t$ is
\begin{eqnarray} q(A) = \begin{pmatrix}a^2 & 2ab & b^2\\ ac &ad+bc & bd
\\ c^2& 2cd &d^2 \end{pmatrix}, \label{qA} \end{eqnarray} whose
eigenvalues are $\det A, \lambda_1^2, \lambda_2^2$. The vector $(2b,
d-a, -2c)^t$ is an eigenvector for $q(A)$ with eigenvalue $\det A$,
and is invariant by $q(A)$ if $\det A = 1$.

The restriction to $\SL(2,\ZZ)$ of the map $A \to q(A)$ defines an
isomorphism onto a discrete subgroup $\Lambda_0$ of automorphisms of
$\SL(3,\ZZ)$ whose each element is non ergodic (each element $q(A)$
leaves fixed a non trivial character of the torus $\TT^3$), but
which acts ergodically on $\TT^3$, since the orbits of the
transposed action on $\ZZ^3 \setminus \{0\}$ are infinite.

\vskip 3mm{\it Extension to a nilmanifold}

Now we extend the action of $\Lambda_0$ to a
nilmanifold. This group is ergodic by Theorem \ref{main}.
Let us consider the real Heisenberg group $H_{2d+1}$ of
dimension $2d+1$, $d \geq 1$, identified with the group
of matrices $(d+2) \times (d+2)$ of the
form:
$$\begin{pmatrix} 1 & x & z \\ 0 & I_d & y \\ 0 & 0 & 1 \end{pmatrix},$$
where $x$ and $y$ are respectively line and column vectors
of dimension $d$, $z$ a scalar and $I_d$ the identity matrix of
dimension $d$. The law of group in $H_{2d+1}$ can be defined by:
$$(x,y,z).(x',y',z') = (x+x',y+y',z+z'+ \langle x,y'\rangle -
\langle x',y\rangle).$$
The map $(x,y,z) \rightarrow (Dx,{}^tD^{-1}y, z)$, for
$D \in \SL(d,\RR)$, defines a group of automorphisms of $H_{2d+1}$.
If the matrices are in $\SL(d,\ZZ)$, these automorphisms preserve the
subgroup $D_{2d+1}$ of elements of $H_{2d+1}$ with integral
coefficients.

The group $\{q(A), A \in \SL(2,\ZZ)\}$ defined by (\ref{qA}) defines
a group $\Gamma$ of automorphisms of the nilmanifold $N/D =
H_{7} /D_7$, with $\Gamma = \{ \tau_A, A \in \SL(2, \ZZ)\}$, where
$\tau_A (x,y,z) = (q(A) x,\, {}^tq(A)^{-1}y,\, z)$. The group
$\Gamma$ acts ergodically on the torus quotient $N/N'D = \TT^3
\times \TT^3$, but each automorphism $(x,y) \rightarrow (q(A)x,\,
{}^tq(A)^{-1}y)$ is non ergodic on $\TT^3 \times \TT^3$.

\section{The spectral gap property \label{SGProp}}

\vskip 3mm
We will now describe some classes of examples where property (SG) is
satisfied.

\vskip 3mm {\it Tori}

As a basic example where property (SG) is valid, let us consider as
in \cite{FuSh99} (See also \cite{Gu06}) the $d$-dimensional torus $X = \TT^d$ endowed with
the Lebesgue measure, and the action of $\SL(d, \ZZ)$ on $\TT^d$ by
automorphisms. The Lebesgue measure is preserved. Every $\gamma \in
\SL(d, \ZZ)$ acts by duality on $\ZZ^d$ by $\gamma^t$. We denote by $\mu^t$ the
push-forward of a probability measure $\mu$ on $\SL(d, \ZZ)$ by the
map $\gamma \to \gamma^t$.

\begin{prop} \label{SGTorus} Let $\mu$ be a probability measure on $\SL(d, \ZZ)$
such that $\supp(\mu^t)$ has no invariant measure on the projective
space $\P^{d-1}$. Let $P$ be the Markov operator on $\TT^d$ defined by
$$P \varphi(x) = \sum_\gamma \varphi(\gamma x) \, \mu(\gamma).$$
Then the corresponding contraction $\Pi_0$ on $L_0^2(X)$ satisfies
$r(\Pi_0) < 1$.
\end{prop} \proof \ The Plancherel formula gives an
isometry $\cal{I}$ between  $L_0^2(\T^d)$ and $\ell^2(\ZZ^d
\setminus \{0 \})$. For $\gamma \in \SL(d, \ZZ)$ we have $\cal{I}
\circ \gamma = \gamma^t \circ \cal{I}$. Hence if $L$ denotes the
convolution operator on $\ell^2(\ZZ^d \setminus \{0 \})$ defined by
$\mu^t$ we have $r(\Pi_0) = r(L)$.

Suppose $r(L) = r(\Pi_0) = 1$ and let $e^{i\theta}$ be a spectral
value. Then two cases can occur. Either there exists a sequence
$(f_n) \in \ell^2(\ZZ^d \setminus \{0 \})$ with $\|f_n\|_2 = 1$ and
$\lim_n \|L f_n - e^{i\theta} f_n\|_2= 0$, or, for some $f \in
\ell^2(\ZZ^d \setminus \{0 \})$ with $\|f\|_2 = 1$, $L^*f = e^{-i\theta} f$.

Since $e^{i\theta} L^*$ is a contraction on $\ell^2(\ZZ^d \setminus
\{0 \})$, its fixed points are also fixed points of its adjoint
$e^{-i\theta} L$, hence $Lf = e^{i\theta} f$. It follows that it
suffices to consider the first case. The condition $\lim_n \| L f_n
- e^{i\theta} f_n\|_2 = 0$ implies $\lim_n \| L f_n\|_2 = \lim_n \|
e^{i\theta} f_n\|_2 = 1$ and also
$$\lim_n [\| L f_n\| _2^2 + \|f_n\|_2
- 2 \Ree \langle L f_n, e^{i\theta } f_n\rangle] =0.$$ Hence $\lim_n
\langle L f_n, e^{i\theta } f_n\rangle = 1$. Since $\langle L f_n,
e^{i\theta} f_n\rangle = \sum_{\gamma} \mu(\gamma) \langle f_n \circ
\gamma^t, e^{i\theta} f_n \rangle$ and $|\langle f_n \circ \gamma^t,
e^{i\theta} f_n\rangle| \leq \|f_n \circ \gamma^t\|_2 \|f_n\|_2$, we
get that, for every $\gamma \in \supp(\mu)$,
$$\lim_n \langle f_n \circ \gamma^t, e^{i\theta} f_n \rangle = 1.$$
Since $|\langle f_n \circ \gamma^t, e^{i\theta} f_n\rangle| \leq
\langle |f_n| \circ \gamma^t, |f_n|\rangle$, we have also $\lim_n
\langle |f_n| \circ \gamma^t, |f_n| \rangle = 1$, hence
$$\lim_n \||f_n| \circ \gamma^t - |f_n| \|_2 = 0.$$
The inequality
$$\||f_n|^2 \circ \gamma^t - |f_n|^2 \|_1 \leq \| |f_n| \circ \gamma^t - |f_n| \|_2
\ \| |f_n| \circ \gamma^t + |f_n| \|_2$$ implies $\lim_n \|
|f_n|^2 \circ \gamma^t - |f_n|^2 \|_1 = 0$.

In other words, if $\nu_n$ denotes the probability measure on
$\ell^2(\ZZ^d \setminus \{0\})$ with density $|f_n|^2$, we have in variational norm:
\begin{equation}
\lim_n \|(\gamma^t)^{-1} \nu_n - \nu_n\| = 0. \label{limvariation}
\end{equation}

Let $\overline \nu_n$ be the projection of $\nu_n$ on $\P^{d-1}$ and $\overline \nu$ a
weak limit of $\nu_n$. By (\ref{limvariation}), we have $(\gamma^t)^{-1}
\overline \nu = \overline \nu$, hence $\gamma^t \overline \nu = \overline \nu,
\forall \gamma^t \in \supp(\mu^t)$, which contradicts the hypothesis on $\supp(\mu^t)$. \eop

\vskip 3mm {\it Remarks} 1) The hypothesis on the support of $\mu^t$ is
satisfied if the group generated by $\supp(\mu)$ has no irreducible
subgroup of finite index. This is a consequence of the following
fact observed by H. Furstenberg: if a linear group has an invariant
measure on $\P^{d-1}$, then either it is bounded or it has a
reducible finite index subgroup (See \cite{Zi84}, p. 39, for a
proof).

\vskip 3mm 2) The above corresponds to a special case
in the characterization of property (SG) given in \cite{BeGu11}, Theorem 5,
for affine maps of $\TT^d$. In particular if $\hat \mu$ is a probability measure
on the group $\Aut(\TT^d) \ltimes \TT^d$ and $\mu$ is its projection in $\Aut(\TT^d)$,
property (SG) for $\hat \mu$ acting on $\TT^d$ is valid if the group generated by $\supp(\mu)$
is non virtually abelian and its action on $\RR^d$ is $ \QQ$-irreducible.

\vskip 3mm {\it Nilmanifolds}

As in Section \ref{nilpSection}, let $N$ be a simply connected nilpotent group, $D$ a lattice in $N$,
$X = N/D$ the corresponding nilmanifold, and $T$ the maximal torus factor.
Let $\mu$ be a probability measure on $\Aut(X) \ltimes N$,  $\overline \mu$
its projection on $\Aut(T) \ltimes T$. Then it is shown in \cite{BeGu11}, Theorem 1,
that if the convolution operator on $L_0^2(T)$ associated with $\overline \mu$ satisfies
property (SG), then the same is true for the operator on $X$ associated with $\mu$.

\vskip 3mm
In view of the torus situation described above, this gives various examples of measures
$\mu$ where property (SG) is valid. If $N$ is a Heisenberg group, more precise
results are available, which will be recalled in Section \ref{examples}.

\vskip 3mm {\it Simple Lie groups}

Let us consider a non compact simple Lie group $G$ and let $\Delta$ be a lattice in $G$,
i.e. a discrete subgroup such that $X = G/\Delta$ has finite volume for the Haar measure
$v$. Let $\mu$ be a probability measure on $G$. It follows from Theorem 6.10
in \cite{FuSh99} that, if $\mu$ is not supported on a coset of a closed amenable
subgroup of $G$, then property (SG) is valid for the action of $\mu$ on $X$.

\vskip 3mm
{\it Compact Lie groups}

We take $X = SU(d)$, $\nu$ the Haar measure on $X$. Then it is known
(See \cite{GaJaSa99}, \cite{BoGa10}) that for $d \geq 3$, if $\Gamma
\subset SU(d)$ is a countable dense subgroup such that the
coefficients of every element of $\Gamma$ are algebraic over $\QQ$
and $\mu$ generates $\Gamma$, then the natural representation of
$\Gamma$ in $L_0^2(X)$ does not contain weakly $\Id_\Gamma$ (cf.
Definitions \ref{def01}).

In particular there are dense free subgroups of $SU(d)$ as above.
Also, if $X =SO(d)$ and $d \geq 5$, there are countable dense
subgroups of $SO(d)$ which have property $(T)$. For example, if
$\Gamma$ is the group of $d \times d$ matrices with coefficients in
$\ZZ[\sqrt 2]$ which preserve the quadratic form $q(x) = \sum_{i=1}^
{d-2} x_i^2 + \sqrt 2\, (x_{d-1}^2 +x_d^2)$, $\Gamma$ has property
$(T)$ (see \cite{Ma91}, p. 136, for similar examples).

Let $A$ be a finite set of generators for $\Gamma$ and $\mu$ a
probability with $\supp(\mu) = A$. Then property (SG) is valid for
the convolution action of $\mu$ on $X$, since $\Gamma$ is ergodic on
$(X, \nu)$, a consequence of the density of $\Gamma$.

\vskip 4mm
\section{Applications of the spectral gap property}\label{appli}
\subsection{Extensions of group actions and random walks}\label{extension}

As in the introduction, let $X$ be a metric space, $\Gamma$ a
countable group of invertible Borel maps of $X$ into itself which
preserve a probability measure $\nu$ on $X$, and $\mu$ a probability
measure on $\Gamma$ with finite support such that $A:= \supp(\mu)$
generates $\Gamma$ as a group. We assume that the action of $\Gamma$
on $(X, \nu)$ is ergodic. We will use both notations: $\sum_a f(ax)
\, \mu(a)$ or $\int f(ax) \,d\mu(a)$.

\vskip 3mm
Let us consider the product space $\Omega = A^{\NN^*}$, with $\NN^*
= \NN \setminus \{0\}$, the shift $\sigma$ on $\Omega$ and the
product measure $\PP = \mu^{\otimes \NN^*}$ on $\Omega$. For $\omega \in \Omega$
we write $\omega =(a_1(\omega), a_2(\omega), ...)$. The
extended shift $\sigma_1$ is defined on $Y= \Omega \times X$  by
$\sigma_1(\omega, x) = (\sigma \omega, a_1(\omega) x)$. Clearly $\sigma_1$
preserves the measure $\PP_1 = \PP \otimes \nu$.

We consider also the bilateral shift on $\hat \Omega := A^\ZZ$ still
denoted by $\sigma$. It preserves the product measure $\hat \PP = \mu^{\otimes \ZZ}$.

\begin{lem} The system $(Y, \PP_1, \sigma_1)$ is ergodic.
\end{lem} \proof \ The dual operator of the composition by $\sigma_1$ on $L^2(\PP_1)$ is
$$Rg(\omega, x) = \sum_{b \in A} g(b\omega, b^{-1} x) \, \mu(b).$$
On functions of the form $(f\otimes \varphi)(\omega, x) = f(\omega) \, \varphi(x)$,
the iterates of $R$ reads
$$R^n (f \otimes \varphi)(x) = \sum_{b_1, ..., b_n \in A^n} f(b_1...b_n \, \omega) \,
\varphi(b_1^{-1}...b_n^{-1} \, x) \, \mu(b_1) ... \mu(b_n).$$

To prove the ergodicity of $\sigma_1$, it suffices to test the convergence of the means
$$\lim_N {1\over N} \sum_{k= 0}^{N-1} R^k g = \int \int g \, d\PP d\nu,$$ when
$g$ is of the form $g(\omega, x) = f(\omega) \varphi(x)$, where
$\varphi$ is in $L^\infty(X)$ and $f$ on $\Omega$ depends only on the first
$p$ coordinates, for some $p \geq 0$. Setting
$$F_{f,\varphi}(x) = \sum_{b_1, ..., b_p \in A^p}
f(b_1...b_p) \, \varphi(b_1^{-1}...b_p^{-1} x) \, \mu(b_1)... \mu(b_p),$$
we have, for $k \geq p$, $R^k (f\otimes \varphi)(\omega, x) = \check P^{k-p}F_{f,\varphi}(x)$,
where $\check P$ is defined by
$$\check P \psi(x) = \sum_{b \in A} \psi(b^{-1} x) \, \mu(b).$$
Ergodicity of $P$, hence of $\check P$, implies the convergence of the means
${1\over N} \sum_{k= 0}^{N-1} R^k (f\otimes \varphi)(\omega, x)$ to the constant
$\int F_{f,\varphi} \, \nu(x) = (\int f \, d\PP(\omega)) \, (\int \varphi \, d\nu(x)).$
\eop

\vskip 3mm
\goodbreak {\it Displacement}

Let $V = \RR^d$ be the $d$-dimensional euclidian space ($d \geq 1$)
and let be given, for each $a \in A$, a bounded Borel map
$x \to c_a(x) = c(a, x)$ from $X$ to $V$. We will call $(c_a(x), a \in A)$
a "displacement".

\vskip  3mm The {\it centering condition} of the displacement is assumed, i.e.
\begin{eqnarray}
\int (\sum_a c_a(x) \,\mu(a)) \, d\nu(x) = 0. \label{centered}
\end{eqnarray}

For $a \in A$ and $(x,v) \in X \times V$, we write $\tilde a(x, v) =
(a x, v + c_a(x))$. Then $\tilde a$ defines an invertible map from
$X \times V$ into itself with ${(\tilde a)}^{-1} (x,v) = (a^{-1} x,
v - c_a(a^{-1} x))$. We denote by $\tilde \Gamma$ the group of Borel
maps of $X \times V$ generated by $\tilde A = \{\tilde a, a \in
A\}$. The action of $\tilde \Gamma$ preserves the fibering of $X
\times V$ over $X$, and the projection of $X \times V$ on $X$ is
equivariant with respect to the action of $\Gamma$ on $X$. We have a
homomorphism $\tilde \gamma \to \gamma$ from $\tilde \Gamma$ to
$\Gamma$ which maps $\tilde a$ to $a$, for every $a \in A$.

In other words, using the displacement $(c_a(x), a \in A)$,
we can extend the action of the group $\Gamma$ on $X$ to the action of the group
$\tilde \Gamma$ generated by the maps $\tilde a$, $a \in A$,
on $X \times V$. Clearly the maps $\tilde \gamma \in \tilde \Gamma$ commute with the translations
on the second coordinate on $X \times V$ by elements of $V$ and therefore are of the form
$\tilde \gamma (x, v) = (\gamma x, v + c(\tilde \gamma, x))$, where $c(\tilde \gamma, x)$ is a map from
$\tilde \Gamma \times X$ to $V$ which satisfies the relation
$$c(\tilde \gamma \tilde \gamma',x) = c(\tilde \gamma, \tilde \gamma' x)
+ c(\tilde \gamma',x), \ \forall \tilde \gamma, \tilde \gamma' \in \tilde \Gamma.$$

For $\tilde \gamma = \tilde a_r...\tilde a_1$, $r \in \NN^*$, we have
$\tilde \gamma  (x, v) =  (a_r...a_1 x, v +
c(\tilde \gamma, x))$, with
\begin{eqnarray}
&c(\tilde \gamma, x) = c(a_1,x)+ c(a_2,a_1x)+... +
c(a_r,a_{r-1}...a_2a_1x). \label{semigroupcocy}
\end{eqnarray}

The displacement satisfies the {\it cocycle property} (for $\Gamma$) if the value
of the sum in (\ref{semigroupcocy}) depends only on the value of the
product $\gamma = a_r...a_1$ in $\Gamma$.

It should be noticed that this cocycle property in general is not
satisfied, since the value of the sum in (\ref{semigroupcocy})
depends in the general case on the "path" $(a_1,...,a_r)$. This is the case
in particular if there is $a$ with
$a, a^{-1} \in A$ and ${\tilde a}^{-1}\not = \widetilde{a^{-1}}$.

A special case is when $(c_a(x), a \in A)$ is a {\it coboundary},
i.e. when there exists $d(x)$ measurable such that $c_a(x) = d(ax) -
d(x), \forall a \in A$. The cocycle property then trivially holds in
$\Gamma$. This is also the case if $c_a(x)$ is a limit of
coboundaries.

\vskip 3mm
{\it Extension of the random walk}

We consider the random walk on $E := X \times V$ defined by the
probability measure $\mu$ and the maps $\tilde a$. Its Markov
operator $\tilde P$ extends the Markov operator $P$ of the random
walk on $X$ given by (\ref{MarkovP}) and is defined by
\begin{eqnarray}
(\tilde P\psi)(x,v) =  \sum_{a \in A} \psi(\tilde a(x,v)) \, \mu(a)
= \sum_{a \in A} \psi(a x, v + c_a(x)) \, \mu(a).
\label{MarkovPTilde}
\end{eqnarray}

Such Markov chains have been considered in the literature under
various names: random walk with internal degree of freedom if $X$ is
finite (\cite{KrSz84}), covering Markov chain (\cite{Ka95}), Markov
additive process (\cite{Uc07}), etc. Intuitively the random walker
moves on $V$ with possible jumps $c_a(x)$, $a \in A$, where $x$
represents the memory of the random walker. Here the steps are
chosen according to the probability $\mu(a)$ which depends on $a$
only. A more general scheme would be to choose the steps $c_a(x)$
according to a weight depending on $(x,a)$. Under spectral gap
conditions on certain functional spaces, it is possible to develop a
detailed study of the iteration $\tilde P^n$ of $\tilde P$ (See for
example \cite{Gu06} when the functional spaces are Hölder spaces).
In the framework of the present paper, no regularity is assumed. We
supposed only that the displacement consists in bounded Borel maps.

\vskip 3mm Recall that $Y =\Omega \times X$. Writing $y = (\omega,
x)$, we define the extension $\tilde \sigma$ of $\sigma_1$ on $Y
\times V$ by $\tilde \sigma (y,v) = (\sigma_1 y, v + c(a_1,x))$. It
preserves the measure $\PP_1 \otimes \ell$, where $\ell$ denotes the
Lebesgue measure on $V$.

The set $Y$ (resp. $Y \times V$) can be identified with the space of
trajectories of the Markov chain defined by $P$ (resp. $\tilde P$).
With the notation of (\ref{semigroupcocy}) we have
$$S_n(y) = S_n(\omega,x) = \sum_{k= 1}^n c(a_k(\omega),
a_{k-1}(\omega)... a_1(\omega) x).$$

Hence, $S_n(y)$ appears as a Birkhoff sum over $(Y, \sigma_1)$.
The iterates of $\tilde \sigma$ on $Y \times V$ read:
$$\tilde \sigma^n(y,v) = (\sigma_1^n y, v + S_n(y)), \forall n \geq 1.$$

Also if we denote by $\tilde \mu$ the push-forward of $\mu$ by the map
$a \to \tilde a$, we can express the iterate $\tilde P^n$ of $\tilde P$ as
$$\tilde P^n \psi(x,v) = \int \psi(\tilde \gamma(x,v)) \ d \tilde \mu^n(\tilde \gamma),$$
where $\tilde \mu^n$ is the $n$-fold convolution product of
$\tilde \mu$ by itself.

Here we are interested in the asymptotic properties of $\tilde
\sigma^n$ and $(S_n(y))$ with respect to the measures $\PP_1 \otimes
\ell$ and $\PP_1$ under the condition that $P$ has "nice" spectral
properties on $X$ (see below). The $L^2$-spectral gap condition can
be compared to the so-called Doeblin condition for the Markov
operator $P$.

\vskip 3mm
The natural invertible extension of $(\Omega \times X  \times V, \tilde \sigma, \PP_1
\otimes \ell)$ is $(\hat \Omega \times X \times V, \tilde \sigma,
\hat \PP_1 \otimes \ell)$, with as above $\tilde \sigma(\omega, x,
v) = (\sigma \omega, a_1(\omega) x, v + c(a_1(\omega),x))$, and where
$\sigma$ is now acting on the bilateral space $\hat \Omega$.

\vskip 3mm
We will need to analyze ${\tilde P}^n$ using methods of Fourier analysis.
Hence we are led to introduce a family of operators $P_\lambda$ on $L^2(X)$,
$\lambda \in V$, defined by
\begin{eqnarray}
P_\lambda \varphi(x) = \sum_{a \in A} e^{i \langle \lambda, \, c_a(x)\rangle} \,
\varphi(ax) \, \mu(a). \label{Plambda}
\end{eqnarray}

We observe that, since $\sup_{a \in A} \|c_a\|_\infty= c < \infty$,
the above formula still makes sense if $\lambda \in \RR^d$ is replaced by
$z \in \CC^d$, and we obtain an operator valued holomorphic function
$z \to P_z$ satisfying, for any $\varphi, \psi \in L^2(X)$,
$$|\langle P_z\varphi, \psi\rangle| \leq e^{c \, \Re e z}
\langle P |\varphi|, |\psi| \rangle \leq e^{c \, \Re e z} \|\varphi\|_2 \|\psi\|_2.$$

This will allow us to use perturbation theory (See \cite{GuHa88} for an analogous situation).

\vskip 3mm The $V$-valued function $h(x) := \sum_{a \in A}
c_a(x) \mu(a)$ belongs to $L_0^2(X)$, in view of the centering condition (\ref{centered}).

Since $r(\Pi_0) < 1$, the restriction of $P-I$ to $L_0^2(X)$ is
invertible, hence we can solve the equation $(P- I)u = h$, with $u
\in L_0^2(X)$. The modified displacement $c'(a, x)$ defined, for $a
\in \Gamma$, by
$$c'(a, x) := c(a, x) - (u(a x) -u(x))$$ satisfies $\nu$-a.e.
\begin{equation}
\sum_{a \in A} c'(a, x) \mu(a) = 0. \label{mart1}
\end{equation}

\vskip 3mm A basic tool for the study of $\tilde P^n$ will be the
analysis of the Fourier operators $P_\lambda$, and in fact their
spectral gap properties. Their family satisfies (as in \cite{GuHa88}
(Lemmas 1 and 2)):
\begin{lem} \label{perturbDev} For any $\varphi \in L^2(X)$, we have
$$P_\lambda \varphi(x) = P\varphi(x) + i \int \langle \lambda,
c_a(x)\rangle \, \varphi(a x) \, d\mu(a) - {1\over 2} \int \langle
\lambda, c_a(x)\rangle^2 \, \varphi(a x) \, d\mu(a)  + |\lambda|^2
o(\lambda).$$

For $\lambda$  small, $P_\lambda$ has a dominant eigenvalue
$k(\lambda)$ given by
$$k(\lambda) = 1 - {1\over 2} \Sigma(\lambda) + |\lambda|^2
o(\lambda),$$ where $$\Sigma(\lambda) = \int \int \langle \lambda,
c'_a(x)\rangle^2 \, d\mu(a) d\nu(x).$$
\end{lem}

In order to analyze more closely the operators $P_\lambda$, we
introduce some definitions related to the aperiodicity of the
displacement.
\begin{defns} \label{33} {\rm We say that the displacement $(c_a, a \in A)$
satisfies \vskip 2mm - (NR) (non reducibility): if there does not
exist $\lambda \in V$, $\lambda \not = 0$, and $d \in L^2(X, \RR^d)$
such that $\nu$-a.e.
$$\langle \lambda, c(a,x) \rangle = \langle \lambda, d(a x) - d(x) \rangle,
\ \forall a \in A ;$$ - (AP) (aperiodicity): if there does not exist
$(\lambda, \theta) \in V \times \RR$, $\lambda \not = 0$, and $d(x)$
measurable, with $|d(x)| = 1$ such that $\nu$-a.e.
\begin{eqnarray*}
e^{i \langle \lambda, c(a,x)\rangle} &=& e^{i\theta} d(a x) / d(x),
\forall a \in A.
\end{eqnarray*}
}\end{defns}

We observe that the formula $\tilde a(x,z) = (ax, z \, e^{i (\langle
\lambda, \, c(a,x) \rangle - \theta)})$ defines an action of $\tilde
\Gamma$ on $X \times \TT$ and that ergodicity of these actions for
every $\lambda \not = 0$ implies condition (AP). Clearly (AP)
implies (NR).

\vskip 3mm
There are special cases (corresponding to functional of Markov chains)
where the previous conditions can be easily verified.
\begin{lem}\label{cnotdep} Assume  that $c_a(x) = c(x)$ does not
depend on $a \in A$, and that the $\RR^d$-valued function $c$ is bounded, not $\nu$-a.e. 0
and satisfies the centering condition $\int c(x) \ d\nu(x) = 0$.

1) If $A \subset \Gamma$ has a symmetric subset $B$ such that $B^2$
acts ergodically on $(X, \nu)$, then (NR) is satisfied. \hfill
\break 2) If the measure $c(\nu)$ is not supported on a coset of a
proper closed subgroup of $\RR^d$, then $(\AP)$ is satisfied.
\end{lem} \proof \
1) If there exists a $\RR^d$-valued functions $d(x)$ such that
$\nu$-a.e. $c(x) = d(ax) - d(x)$, for any $a \in A$, then, for any
$a,a' \in A$, $d(ax) = d(a'x)$, i.e. $d(\gamma x) = d(x)$ for any
$\gamma \in A A^{-1}$. Since the subgroup generated by $A A^{-1}$ is
ergodic on $(X, \nu)$ we get that, for some $c_0 \in \RR^d$, $d(x)=
c_0$, hence $c(x) = 0$, $\nu$-a.e.

\vskip 3mm 2) If $(\AP)$ does not hold, there exist $\lambda \not =
0$ in $\RR^d$, $\theta_\lambda \in \RR$ and a cocycle
$\sigma_\lambda(\gamma,x)$ on $\Gamma \times X$ with values in the
group of complex of modulus 1, such that for $\nu$-a.e. $x$
$$\sigma_\lambda(a,x) = e^{i(\langle \lambda, c(x)\rangle -
\theta_\lambda)}, \forall a \in A.$$

In particular, taking $a,a^{-1} \in B$, we have
$$1 = \sigma_\lambda(a^{-1},ax) \, \sigma_\lambda(a,x) =
e^{i(\langle \lambda, c(x)+c(ax)\rangle - 2\theta_\lambda)};$$
hence, for any $a,a' \in B$, $\nu$-a.e. $e^{i\langle \lambda,
c(ax)\rangle} = e^{i\langle \lambda, c(a'x) \rangle}$.

Since $B^2$ acts ergodically on $(X, \nu)$ we get, for  $\lambda
\not = 0$ and some $c_\lambda$ of modulus 1, $e^{i(\langle \lambda,
c(x)\rangle} = c_\lambda$. This means that $c(x)$ belongs to the
coset of the proper closed subgroup of $\RR^d$ defined by
$e^{i(\langle \lambda, v \rangle} = c_\lambda$, which contradicts
the hypothesis. \eop

\subsection{Ergodicity, recurrence/transience}

In this section, we study ergodicity, recurrence and transience of
the extended dynamical systems considered above. First we recall
briefly the notion of recurrence in the framework of dynamical systems.
(See (\cite{Sc77}, \cite{Aa07}.)

Let $(Y,\lambda, \tau)$ be a dynamical system with $Y$ a metric space,
$\lambda$ a probability measure on $Y$ and $\tau$ an invertible Borel map of $Y$
into itself which preserves $\lambda$. We suppose the system ergodic.
If $\varphi$ is a Borel map from $Y$ to $\RR^d$,
the ergodic sums $S_n \varphi (y) = \sum_{k=0}^{n-1}\varphi(\tau^k y)$
define a "random walk in $\RR^d$ over the dynamical system" $(Y,\lambda, \tau)$.
The corresponding skew product is the dynamical system defined on
$(Y \times \RR^d, \lambda \times \ell)$ by the transformation
$\tau_\varphi : (y, v) \to ( \tau y, v + \varphi(y))$.

We say that $y \in Y$ is {\it recurrent} if,
for every neighborhood $U$ of $0$ in $\RR^d$,
$$\sum_{n \geq 0} 1_U(S_n \varphi(y)) = +\infty.$$
We say that $y$ is {\it transient} if, for every neighborhood $U$ of
$0$ in $\RR^d$ this sum is finite. The cocycle $(S_n \varphi)$ is
{\it recurrent} if a.e. point $y \in Y$ is recurrent. It is
transient if a.e. point $y \in Y$ is transient.

Since the set of recurrent points is clearly invariant
and the system $(Y, \lambda, \tau)$ is ergodic,
every cocycle $(S_n\varphi)$ is  either transient or recurrent.

For the sake of completeness, let us give a simple proof
of the following known equivalence:

\vskip 3mm \begin{proposition} The recurrence of
$(S_n \varphi)$ is equivalent to the
conservativity of the system $(Y \times \RR^d, \lambda \otimes \ell,
\tau_\varphi)$. \end{proposition} \proof \ By definition the dynamical system
$(Y \times \RR^d, \lambda \otimes \ell, \tau_\varphi)$ is conservative if, for every
measurable set $A$ in $Y \times \RR^d$ with positive measure, for a.e. $(y,v)\in
A$ there exists $n \geq 1$ such
that $\tau_\varphi^n(y,v) = (\tau^n y, v + S_n\varphi(y)) \in A$.

We show conversely that this property holds if
$(S_n \varphi)_{n \geq 1}$ is recurrent in the sense of the above definition.

We can suppose that $A$ is included in $Y \times L$,
where $L$ is a compact set in $\RR^d$. One checks easily that
the set $B:=\{(y,v) \in A: \forall n \ge 1,
\tau_\varphi^n(y,v) \not \in A\}$ has pairwise disjoint images.

Using the recurrence of $(S_n\varphi)$, one can find for every $\varepsilon >0$
a compact set $K_\varepsilon$ such that, for a set of measure $\geq 1 - \varepsilon$
of points $y$, the sums $S_n\varphi(y)$ belongs to $K_\varepsilon$ for infinitely many
$n$ (we use the fact that for a.e. $y$, the set
$\{S_n\varphi(y), n \geq 1\}$ has a finite accumulation point and that there is a neighborhood
of this accumulation point in which $S_n\varphi(y)$ returns infinitely often).

The measure of the set $F_\varepsilon:=Y \times (L+K_\varepsilon)$ is finite.
Since $B$ has pairwise disjoint images, we have
\begin{eqnarray*}
&&\int_Y 1_B(y,v) \sum_{n \geq 1} 1_{F_\varepsilon}(\tau_{\varphi}^n (y,v)) \ \lambda(dy) \, d\ell(v) =
\sum_{n \geq 1} (\lambda \times \ell)(B \cap \tau_{\varphi}^{-n} F_\varepsilon) \\
&=& \sum_{n \geq 1} (\lambda \times \ell)(\tau_{\varphi}^n B \cap F_\varepsilon)
\leq (\lambda \times \ell)(F_\varepsilon) < +\infty,
\end{eqnarray*}
and therefore, for a.e.  $(y,v)$ in $B$, $\sum_{n \geq 1}
1_{F_\varepsilon}(\tau_{\varphi}^n (y,v)) <\infty$.

This implies that $(\lambda \times \ell)(B) \leq \varepsilon$, hence $B$ has
measure 0. The set $A$ satisfies the announced property. \eop

\vskip 3mm
From the proposition it follows that, when $(S_n \varphi)$ is recurrent,
the random walk visits $\PP\otimes \nu$-a.e. infinitely often any
neighborhood of 0 in $V$, i.e., $\liminf_{n\to \infty}
\|S_n(\omega, x)\| = 0$, $\PP\otimes \nu$-a.e.
When the random walk $(S_n \varphi)$ is transient, then  $\lim_{n\to \infty}
\|S_n(\omega, x)\| = \infty$, $\PP\otimes \nu$-a.e. on $\Omega
\times X$.

In the transient case, the system $(Y \times \RR^d, \lambda \otimes \ell,
\tau_\varphi)$ is dissipative, i.e. there exists a Borel subset
$B \subset Y \times \RR^d,$ such that
$$Y \times \RR^d = \bigcup_{n \in \ZZ} \tau_\varphi^n B {\rm \ and \ }
(\lambda \otimes \ell) (\tau_\varphi^n B \cap B) = 0, \ \forall n \in \ZZ \setminus
\{0\}.$$

Now we will study these properties of recurrence and transience in the case of the
random walk $(S_n(\omega, x))$ over the dynamical system
$(\hat \Omega \times X, \hat \PP_1, \sigma_1)$ and its extension
$(\hat \Omega \times  E, \tilde \sigma, \hat \PP_1 \otimes \ell)$
defined at the beginning of this section.

\begin {thm} \label{recurr} Assume that $(X, \nu, \Gamma, \mu)$
satisfies property (SG).

1a) If the displacement $(c_a, a \in A)$ satisfies (NR), then
$({1 \over \sqrt n} S_n(\omega, x))_{n \geq 1}$ converges in law with
respect to $\PP \otimes \nu$ to the centered  normal law on $V$
with non degenerate covariance $\Sigma$.

1b) If the displacement $(c_a, a \in A)$ satisfies $(\AP)$, then the
local limit theorem holds: for any $\varphi \in L^2(X)$ and $f$
continuous with compact support, $\tilde \alpha = \alpha \nu \otimes
\delta_0$, with $\alpha \in L^2(X)$, $\alpha \geq 0$, $\nu(\alpha) =
1$, we have
\begin{eqnarray}
\lim_n (2 \pi n)^{d/2} (det \,\Sigma)^{1\over 2} \tilde P^n \tilde
\alpha (\varphi\otimes f) =  \nu(\varphi)\, \ell(f). \label{thmloc}
\end{eqnarray}

2a) For $d \leq 2$ $(S_n)$ is recurrent:
$\PP\otimes \nu$ a.e. $\liminf_{n\to \infty} \|S_n(\omega, x)\| = 0$.

2b) For $d \leq 2$, if the displacement $(c_a, a \in A)$ satisfies
$(\AP)$, then $\,\tilde \sigma$ is ergodic with respect to $\hat
\PP_1 \otimes \ell$.

2c) If $d \geq 3$, if the displacement $(c_a, a \in A)$ satisfies
$(\AP)$, $(S_n)$ is transient: $\PP\otimes \nu$ a.e. on $\Omega
\times X$ $\lim_{n\to \infty} \|S_n(\omega, x)\| = \infty$.

3) For any $d \geq 1$, if the displacement $(c_a, a \in A)$
satisfies $(\AP)$, the equation $\tilde P h = h$, $h \in
L^\infty(\nu \otimes \ell)$, has only constant solutions.

\end{thm} \proof \ 1) We have
$$S_n(\omega, x) = \sum_{k=1}^n Y_k(\omega, x) + u( a_n... a_1 x) - u(x),$$
with $Y_k(\omega, x) = c'(a_k, a_{k-1}... a_1 x)$.
We observe that
$$E(Y_k(\omega, x) | a_1, ..., a_{k-1}) = \int c'(a,x) \ d\mu(a) = 0,$$
since $h(x) = \int c(a, x) \, d\mu(a) = (Pu-u)(x)$.

\vskip 3mm On the other hand, for any $v \in V$, by definition of
$\Sigma(v)$ and from (\ref{mart1}) we have the martingale property
and in particular $\EE(\langle Y_k, v \rangle\, \langle Y_\ell, v
\rangle) = 0, {\rm \ if \ } k \not = \ell$, and
\begin{eqnarray*}
\EE(\langle Y_k, v \rangle^2) = \int \,  \langle Y_k(\omega, x), v\rangle^2
d\PP(\omega) \, d\nu(x) = \int \int \langle \lambda,
c'_a(x)\rangle^2 \, d\mu(a) d\nu(x) = \Sigma(v).\\
\end{eqnarray*}
Clearly $Y_k = Y_1 \circ \sigma_1^k$ and $\sigma_1$ preserves the
measure $\PP \otimes \nu$. Ergodicity of $\sigma_1$ and the ergodic theorem imply $\PP
\otimes \nu$-a.e.
$$\lim_{n \to \infty} {1\over n} \sum_{k=1}^n \langle Y_k, v\rangle
^2 = \Sigma(v).$$

Hence
$$\lim_n {1\over n} \EE(\langle \sum_1^n Y_k, v\rangle^2) = \Sigma(v) =
\lim_n {1\over n} \EE(\langle S_n,v\rangle^2).$$

Brown's central limit theorem (\cite{Br71}) applies to $(Y_k)$, and gives the CLT for $(S_n)$
since the coboundary term $u( a_n... a_1 x) - u(x)$ is bounded.
So we get the convergence in law of ${1 \over \sqrt n} S_n(\omega, x)$ with
respect to $\PP \otimes \nu$ to the centered  normal
law on $V$ with covariance $\Sigma$. The non degeneracy of $\Sigma$ follows from the formula
$\Sigma(v) = \int \langle v, c_a'(x) \rangle^2 \, d\mu(a) \, d\nu(x)$ and Condition (NR).

The statement 1b) (the convergence (\ref{thmloc}))
is proved in Lemma \ref{convRadon} below.

\vskip 3mm 2a) Using the CLT as a recurrence criterion for the
$\RR^2$-valued $\ZZ$-cocycle, $(S_ny)$ over the measure preserving
transformation $\sigma_1$ (cf. \cite{Sc98} or \cite{Co99}), the
recurrence property follows:
\begin{equation}
\liminf_{n \to \infty} \|S_n(\omega, x)\| = 0. \label{recTo0}
\end{equation}

2b) We observe that the trajectories of the random walk on
$E = X \times V$ defined by $\mu$ are given by
$$X_n(\omega, x, v) = (\sigma_1^n(\omega, x), v+ S_n(\omega, x)).$$

By (\ref{recTo0}), for any relatively compact open set $U \subset V$, $\PP \otimes
\nu$-a.e. on $\Omega \times X \times U$ there exists $n(\omega, x,
v) \geq 1$ with $X_n(\omega, x, v) \in \Omega \times X \times U$.

In other words, the Markov kernel $\tilde P$ on $E$ satisfies
Property R defined in \cite{GuRa09}. Hence, using Proposition 2.6 in
\cite{GuRa09}, the ergodicity of $(\hat \Omega \times E, \tilde
\sigma, \hat \PP_1 \otimes \ell)$ will follow if we show that the
equation $\tilde P h = h$, for $h \in L^\infty(\nu \otimes \ell)$,
has only constant solutions.

Since $\tilde P h = h$, we have for any $n \in \NN$, $\varphi \in L^2(X)$,
$f \in L^1(V)$ with $\int f(v) \, d\ell(v) = 0$,
$$\langle ({{\tilde P}^*})^n (\varphi \otimes f), h\rangle = \langle\varphi
\otimes f, h\rangle.$$

Lemma \ref{spectralGapProd} below gives $\langle \varphi \otimes
f,h\rangle = 0$, hence $h$ is invariant by translation by $v$ and
defines an element $\overline h \in L_\infty(X, \nu)$ with $P
\overline h = \overline h$. Then we have $\sum_{a\in A} \overline
h(a x) \mu(a) = \overline h(x)$, hence the invariance of $h$ by $a \in A$.
Since $\nu$ is $\Gamma$ ergodic, $\overline h$ is constant
$\nu$-a.e. Therefore $h$ is constant $\nu\otimes \ell$-a.e.
This proves 2b).

\vskip 3mm 2c) Let us show that, for any relatively compact
subset $U$ of $E$, for a.e. $(x,v) \in U$ we have on $U$: $\sum_{n=
1}^\infty 1_U( a_n... a_1 x, v+S_n(\omega, x))) < +\infty$.

We have, for every non negative Borel function $\psi$ on $E$:
$$\EE(\sum_{n= 1}^\infty \psi( a_n... a_1 x, v+S_n(\omega, x))) = \sum_{n= 1}^\infty
{\tilde P}^n\psi(x,v).$$

Here we will prove the convergence $\sum_{n= 1}^\infty \langle
|\tilde P^n \psi, \psi\rangle| < \infty$ for $\psi$ of the form
$\varphi\otimes f$. Since we can choose $\psi \geq 1_U$, this will
implies
$$\EE(\sum_1^\infty 1_U(a_n...a_1 x, v + S_n(\omega, x) 1_U(x,v)) < \infty,$$
hence the result. This convergence follows from Lemma
\ref{convRadon} below.

To finish the proof of 2), we observe that if $\varphi$ is a
continuous function with compact support on $V$ and $\psi = \alpha
\otimes \varphi$, where $\alpha \in L^2(X)$ and $\tilde \alpha =
\alpha \nu \otimes \delta_0$, then
$$\langle \tilde P^n \psi, \psi\rangle = (\tilde P^n \tilde \alpha) (\psi)
\ell(\varphi).$$

In particular by (\ref{thmloc}) the sequence $(n^{d/2} \langle P^n
\psi, \psi \rangle)$ is bounded. If $d > 2$ the series
$\sum_{n=0}^\infty  |\langle \tilde P^n \psi, \psi\rangle|$ converges.
Hence the result.

3) The assertion is shown in the proof of 2b). \eop

Now, under the assumption $(\AP)$ as in the theorem, we prove the
lemmas used in the previous proof.

\begin{lem} \label{spectralGapProd}
For any $\varphi \in L^2(X)$ and any $f \in L^1(V)$
with $\int f(v) \, d\ell(v) = 0$, we have
$$\lim_{n \to \infty} \|{\tilde P}^{*^n}(\varphi \otimes f) \|_1 = 0.$$
\end{lem} \Proof \ In the proof of Proposition 3.6 of \cite{GuSt04}, a Markov
operator $Q$ on $X \times \RR^d$ which commutes with the
$\RR^d$-translations is considered and it is proved that $\lim_{n
\to \infty} \|Q^{n}(u \otimes f) \|_1 = 0$ for Hölder continuous
functions $u$ in $H_\varepsilon(X)$ and $f$ as above. The essential
points of the proof are polynomial growth of $\RR^d$ as a group and
a spectral gap property for the $Q$-action on function of the form
$u \otimes \lambda$ where $u \in H_\varepsilon(X)$ and $\lambda$ is
a character of $\RR^d$.

Here we observe that the adjoint operator ${\tilde P}^*$ of $\tilde P$
on $L_0^2(E, \nu \otimes \ell)$ is associated with $\check \mu$ the symmetric of
$\mu$ which has the same properties as $\mu$ as it was observed
above.

The action of ${\tilde P}^*$ is also well defined on the functions
of the form $u\otimes \lambda$ where $\lambda$ is a fixed  character
of $V$ and $u$ is in $L^2(X)$. It reduces there to the action of
$P_\lambda^*$ on $L^2(X)$ hence using $b) \Rightarrow a)$ of
Proposition \ref{P0} below we get that $(\AP)$ implies that
$P_\lambda^*$ has a spectral gap. Hence the lemma follows from the
proof of Proposition 3.6 in \cite{GuSt04} with $Q = \tilde P^*$.
\eop

\begin{lem} \label{convRadon} Let $\alpha$ be a probability measure
on $X$ which has a $L^2$-density with respect to $\nu$. Let $\tilde
\alpha$ be the probability measure $\alpha \otimes \delta_0$ on $X\times V$
and let $\tilde \mu_n:= (2 \pi n)^{d/2} (det \,\Sigma)^{1\over 2} \tilde P^n \tilde
\alpha$. Then the sequence of Radon measures $(\tilde \mu_n)$ on $X \times
V$ satisfies $\lim_n \tilde \mu_n(\varphi\otimes f) =  \nu(\varphi) \, \ell(f)$
for any $\varphi \in L^2(X)$ and $f$ continuous with compact support.
\end{lem}
\Proof \ Let $\varphi$ be a function in $L^2(X)$ and $f \in
L^1(V)$ be such that its Fourier transform $\hat f(\lambda) = \int
f(v) \, e^{i\langle \lambda, v \rangle} \, d\ell(v)$ has a compact
support on $V$. Then, by the inversion formula we have:
$$f(v) = (2\pi)^{-d} \int \hat f (\lambda)
\, e^{-i\langle \lambda, v \rangle} \, d\ell(\lambda).$$ As in
\cite{Br68}, p. 225, we test the convergence of $\tilde
\mu_n^\varphi (f)=\tilde \mu_n(\varphi \otimes f)$ using functions
$f$ as above. We apply the method of \cite{GuHa88} for proving the
local limit theorem, giving only the main steps of the proof.
According to $b) \Rightarrow a)$ of Proposition \ref{P0} below, we
observe that the operator $P_\lambda$ considered above satisfies
$r(P_\lambda) < 1$ for $\lambda \not = 0$, in view of Condition
$(\AP)$. Furthermore, by perturbation theory, for $\lambda$ small
enough, $P_\lambda$ has a dominant eigenvalue $k(\lambda)$ and a
corresponding one dimensional projection $p_\lambda$ such that:
\begin{eqnarray*}
P_\lambda &=& k(\lambda) p_\lambda + R_\lambda, \\ R_\lambda
p_\lambda &=& p_\lambda  R_\lambda, \ r(R_\lambda) < |k(\lambda)|, \\
k(\lambda) &=& 1 - {1\over 2} \Sigma(\lambda) + |\lambda|^2
o(|\lambda|).
\end{eqnarray*}

Also $p_\lambda$, $r_\lambda$ depend analytically on $\lambda$.
These facts will allow us to adapt the analogous proof of
\cite{GuHa88}. We write $\tilde P^n \tilde \alpha$ as follows:
\begin{eqnarray*}\tilde P^n \tilde \alpha &=& \int \delta_ {\gamma x}
\otimes \delta_{c(\tilde \gamma,x)} \, d {\tilde \mu}^n(\tilde \gamma) \,
d\alpha(x); \\
P^n \tilde \alpha (\varphi \otimes f)
&=& (2\pi)^{-d} \int \, \varphi( \gamma x) \, e^{-i\langle
\lambda, c(\tilde \gamma,x)\rangle} \, \hat f(\lambda)
\, d\tilde \mu^n(\tilde \gamma) \, d\alpha(x) \, d\ell(\lambda)\\
&=& (2\pi)^{-d} \int_{\RR^d} \alpha(P_{-\lambda}^n \varphi) \, \hat
f(\lambda) \, d\ell(\lambda),
\end{eqnarray*}
hence
\begin{eqnarray*}
\tilde \mu_n(\varphi \otimes f) = (2\pi)^{-d/2} (det \,
\Sigma)^{1\over 2} \, n^{d/2} \int_{\RR^d} \alpha(P_{-\lambda}^n
\varphi) \, \hat f(\lambda) \, d\ell(\lambda).
\end{eqnarray*}

Since $r(P_\lambda) < 1$ for $\lambda \not = 0$, the integration can
be reduced, in the limit, to a small neighborhood $U$ of 0 in $\RR^d$
and it suffices to consider $$I_n:= (2\pi)^{-d/2} (det \,
\Sigma)^{1\over 2} \int_{\sqrt n U} \alpha(P_{-\lambda / \sqrt
n}^n(\varphi) \, \hat f(\lambda / \sqrt n) \ d\ell(\lambda);$$

Using the spectral decomposition of $P_\lambda$ we see that $\lim_{n \to
\infty} I_n = \lim_{n \to \infty} J_n$ with
$$J_n= (2\pi)^{-d/2} (det \, \Sigma)^{1\over 2} \int_{\sqrt n U}
k(\lambda / \sqrt n)^n \alpha(P_{\lambda /\sqrt n}(\varphi)) \, \hat
f(\lambda / \sqrt n) \ d\ell(\lambda).$$

Since $\lim_{n \to \infty} k(\lambda /\sqrt n)^n = e^{-{1\over 2}
\Sigma(\lambda)}$ and $\lim_{\lambda \to 0} p_\lambda(\varphi) =
\nu(\varphi)$, we get
$$\lim_{n \to \infty} J_n = (2\pi)^{-d/2} (det \, \Sigma)^{1\over 2}
\int e^{-{1\over 2} \Sigma(\lambda)} \, \nu(\varphi) \, \hat f(0) \
d\ell(\lambda);$$ hence $\lim_{n \to \infty} J_n = (\nu \otimes
\ell)(\varphi \otimes f)$. \eop

\vskip 5mm
\begin{prop} \label{P0} Let $\Pi_0$ be the restriction of $P$ to $L_0^2(X)$
and let $P_\lambda$, $\lambda \in V$, $\lambda \not = 0$, be defined
on $L^2(X)$ by (\ref{Plambda}) with $c_a(x) \in L^\infty(X)$ for
every $a \in \supp(\mu)$. Assume that $r(\Pi_0) < 1$. Then the
following properties are equivalent: \hfill \break a) the spectral
radius $r(P_\lambda)$ of $P_\lambda$ acting on $L^2(X)$ is 1; \hfill
\break b) the condition $(\AP)$ is not satisfied at $\lambda$: there
exists a real $\theta$ and a measurable function $\alpha$ such that
\begin{eqnarray*}
e^{i(\langle \lambda, \, c_a(x) \rangle - \theta)} = e^{i(\alpha(ax)
- \alpha(x))}, \ \nu-a.e.
\end{eqnarray*}
c) there exists $\theta_\lambda \in \R$ such that $e^{i(\langle
\lambda, \, c_a(x) \rangle - \theta_\lambda)}$ extends as a cocycle
$\sigma_\lambda(\gamma,x)$ on $\Gamma \times X$, with values in the
group of complex numbers of modulus 1, and the representation
$\rho_\lambda$ of $\Gamma$ on $L^2(X)$ contains $\Id_\Gamma$, where
\begin{eqnarray}
(\rho_\lambda(\gamma) \varphi)(x) = \sigma_\lambda(\gamma^{-1},x) \,
\varphi(\gamma^{-1} x).\label{rhoRepr}
\end{eqnarray}
\end{prop}
\proof \ a) $\Rightarrow$ b)

We begin as in the proof of Proposition \ref{SGTorus}.
Assume $r(P_\lambda) = 1$ and let $e^{i\theta}$, ($\theta \in [0, 2
\pi[$) be a spectral value of $P_\lambda$. Then, either the subspace
${\rm Im}(e^{i\theta} - P_\lambda)$ is not dense in $L^2(X)$ or there
exists $\varphi_n \in L^2(X)$, with $\|\varphi_n\| = 1$, such that
\begin{eqnarray*}
\lim_n \|P_\lambda \varphi_n - e^{i\theta} \varphi_n \|_2 = 0.
\end{eqnarray*}

In the first case, there exists $\varphi \in L^2(X)$ with
$e^{i\theta} P_\lambda^* \varphi = \varphi.$ Since $e^{-i\theta}
P_\lambda$ is a contraction of $L^2(X)$, this implies $
P_\lambda\varphi = e^{i\theta} \varphi.$ Hence it suffices to
consider the second case. We have:
\begin{eqnarray*}
0 \leq \|P_\lambda \varphi_n - e^{i\theta} \varphi_n \|_2^2 =
\|P_\lambda \varphi_n\|_2^2 + \|\varphi_n\|_2^2 - 2 \Ree \langle
P_\lambda \varphi_n, e^{i\theta} \varphi_n\rangle \leq 2 - 2 \Ree
\langle P_\lambda \varphi_n, e^{i\theta} \varphi_n\rangle.
\end{eqnarray*}

Then the condition $\lim_n \|P_\lambda \varphi_n - e^{i\theta}
\varphi_n \|_2 = 0$ is equivalent to:
$$\lim_n \Ree \langle P_\lambda \varphi_n, e^{i\theta} \varphi_n\rangle = 1,$$
i.e. to
$$\lim_n \langle P_\lambda \varphi_n, e^{i\theta} \varphi_n\rangle
= 1.$$

We have also, since $\|\varphi_n\|_2 = 1$,
$$0 \leq |\langle P_\lambda \varphi_n,
e^{i\theta} \varphi_n\rangle| \leq  \langle P |\varphi_n|,
|\varphi_n| \rangle \ \leq 1.$$
It follows $\lim_n \langle P |\varphi_n|, |\varphi_n| \rangle = 1$,
i.e.
$$ \lim_n \|P |\varphi_n| -|\varphi_n|\|_2 = 0.$$

We can write $|\varphi_n| = \langle |\varphi_n|, 1\rangle + \psi_n$,
where $\psi_n:= (|\varphi_n| - \langle |\varphi_n|, 1\rangle 1 \in
L_0^2(X)$ and therefore $\lim_n \|(P-I) \psi_n \|_2 = 0$ and
$\langle |\varphi_n|, 1 \rangle \leq \|\varphi_n\|_2 \leq 1$.

Since $r(\Pi_0) < 1$, the restriction of $\Pi_0 -I$ to $L_0^2(X)$ is
invertible. Hence the condition $\lim_n \|(P-I) \psi_n\|_2 = 0$
implies $\lim_n \|\psi_n\|_2 = 0$.

If $c \in \RR_+$ is a limit of a subsequence $\langle
|\varphi_{n_i}|, 1\rangle$ of $\langle |\varphi_{n}|, 1\rangle$, we
get $\lim_i \||\varphi_{n_i}| -c\|_2 = 0$.

Since $\|\varphi_n\|_2 = 1$, we have $c =1$, hence the convergence
$\lim_n \||\varphi_n| -1\|_2 = 0$.

On the other hand, the condition $\lim_n \langle P_\lambda
\varphi_n, e^{i\theta} \varphi_n\rangle = 1$ can be written as
$$\lim_n \int \langle e^{i \langle \lambda, c_a \rangle} \varphi_n \,\circ \, a, e^{i\theta}
\varphi_n \rangle \ d\mu(a) \ = 1,$$ where for each $a \in
\supp(\mu)$
$$|\langle e^{i\langle \lambda, c_a \rangle} \varphi_n \, \circ \, a, e^{i\theta}
\varphi_n\rangle| \leq 1.$$

It follows, for any $a \in \supp(\mu)$:
$$\lim_n \|e^{i \langle \lambda, c_a \rangle} \varphi_n \circ a - e^{i\theta}
\varphi_n\|_2 = 0.$$

We can write $\varphi_n(x) = |\varphi_n(x)|e^{i\alpha_n(x)}$, with
$\alpha_n(x) \in [0, 2\pi [$. Hence:
\begin{eqnarray*}
&&e^{\langle \lambda, c_a \rangle} \varphi_n\circ a - e^{i\theta} \varphi_n \\
&&= e^{\langle \lambda, c_a \rangle +i \alpha_n \circ a} (|\varphi_n \circ a| - 1) -
e^{i\theta} (|\varphi_n| - 1) e^{i\alpha_n} + e^{i \langle \lambda, c_a \rangle +i \alpha_n
\circ a} - e^{i(\alpha_n +\theta)}.
\end{eqnarray*}

Hence
$$\lim_n \|e^{i(-\alpha_n + \alpha_n \circ a + \langle \lambda, c_a \rangle - \theta)} - 1
\|_2 = 0,$$
therefore, for a subsequence $(n_k)$
\begin{eqnarray*}
e^{i(\langle \lambda, c_a(x) \rangle - \theta)} = \lim_{k}
e^{i(\alpha_{n_k}(ax) - \alpha_{n_k} (x))}, \ \nu-a.e.
\end{eqnarray*}

Clearly $\lim_{k} e^{i(\alpha_{n_k}(ax) - \alpha_{n_k} (x))}$ is the
restriction to $A \times X$ of a cocycle $\sigma_\lambda(\gamma,x)$
of $\Gamma \times X$.

On the other hand property (SG) implies strong ergodicity of the
action of $\Gamma$ on $X$; hence proposition 2.3 of \cite{Sc80}
gives the existence of a measurable function $\alpha$ on $X$ such
that $\sigma_\lambda(a, x) = e^{i(\langle \lambda, c_a(x) \rangle -
\theta)} = e^{i(\alpha(ax) - \alpha(x))}$, $\nu$-a.e.

\vskip 3mm b) $\Rightarrow$ c)

With $\varphi(x) = e^{-i \alpha(x)}$, we have by condition b)
$\sigma_\lambda(a,x) = \varphi(x) / \varphi(ax)$ which extends to
$\Gamma \times X$ as a cocycle. By the definition of
$\rho_\lambda(a^{-1})$ (cf. (\ref{rhoRepr})), we have
$$(\rho_\lambda(a^{-1}) \, \varphi)(x) = \varphi(x).$$

Since $A$ generates $\Gamma$ as a group, this means that the
representation $\rho_\lambda$ contains $\Id_\Gamma$.

\vskip 3mm c) $\Rightarrow$ a)

Let $\check \mu$ be the push-forward of $\mu$ by the map $\gamma \to
\gamma^{-1}$. We observe that, by the definition of $\rho_\lambda$,
$\rho_\lambda(\check \mu) = e^{-i\theta} P_\lambda$.

This means that $e^{i\theta}$ is an eigenvalue of $P_\lambda$, so
that $r(P_\lambda) = 1$. \eop

\vskip 3mm {\it Remarks:} 1) In general, it is non trivial to
calculate the set of $\lambda \in \RR^d$ such that $r(P_\lambda) =
1$. However, Corollary \ref{discSubG} is useful for this question.

2) Also we observe that condition c) in the proposition implies that
the action of $\Gamma$ on $X \times \TT$ given by $\gamma(x,t) =
(\gamma x, t \sigma_\lambda(\gamma, x))$ is not ergodic.

\vskip 3mm
\begin{cor} \label{discSubG} Assume property (SG). The set
$R_\mu=\{\lambda \in \RR^d: r(P_\lambda) = 1 \}$ is a closed
subgroup of $\RR^d$. It is discrete if (NR) is valid.
\end{cor} \proof \ \ Assume $\lambda_1, \lambda_2$ satisfy
$r(P_{\lambda_1}) = r(P_{\lambda_2}) = 1$. Then condition b) of the
Proposition gives the existence of $\theta_1, \theta_{2} \in \RR$,
$\varphi^1, \varphi^2$ with $|\varphi^1| = |\varphi^2| = 1$, such
that, for $j = 1, 2$ and $\nu$-a.e.
$$e^{i\langle \lambda_j, c_a(x)\rangle - \theta_{\lambda_j}}
=  \varphi^j(ax) / \varphi^j(x).$$ It follows
$$e^{i\langle \lambda_1- \lambda_2, c_a(x)\rangle - \theta_{\lambda_1} + \theta_{\lambda_2}}
= {\varphi^1(ax)\overline {\varphi^2}(ax) \over \varphi^1(x)
\overline {\varphi^2(x)}},$$ i.e. condition b) is satisfied with
$\lambda_1 -\lambda_2$, $\theta_{\lambda_1} - \theta_{\lambda_2}$,
$\varphi^1(x) \overline {\varphi^2(x)}$. Hence $R_\mu$ is a group.

\vskip 3mm
The definition of $P_\lambda$, $P_{\lambda'}$ gives the following
inequality:
\begin{eqnarray*}
&&|(P_\lambda - P_{\lambda'}) \varphi (x)| \leq \sup_{x,a} |e^{i\langle
\lambda - \lambda',\, c_a(x)\rangle} - 1 | \, P|\varphi| (x),
\end{eqnarray*}
hence, since $|c_a(x)|$ is bounded by a constant $c$,
\begin{eqnarray*}
\|(P_\lambda - P_{\lambda'}) \varphi\|_2 \leq \sup_{a \in A} \|c_a\|_\infty\, |\lambda - \lambda'| \, \|\varphi\|_2
\leq  c |\lambda - \lambda'| \| \varphi\|_2.
\end{eqnarray*}

Therefore, we have $\|P_\lambda - P_{\lambda'}\| \leq  c |\lambda - \lambda'|$,
which implies that, if $\lambda$ is fixed with $r(P_\lambda) < 1$, we have
also $r(P_{\lambda'}) < 1$ for $\lambda'$ sufficiently close to
$\lambda$. Hence $R_\mu$ is closed.

Assume now Condition (NR) is satisfied. Observe that $\|P_\lambda -
P_0\| \leq  c |\lambda|$. Since $r(\Pi_0) < 1$, the spectrum of $P =
P_0$ consists of $\{1 \}$ and a compact subset of the open unit
disk. Hence $P$ has a dominant isolated eigenvalue, which is a
simple eigenvalue. By perturbation theory, this property remains
valid in a neighborhood of 0.

Using Lemma \ref{perturbDev} and the fact (noted in Theorem \ref{recurr} part 1) that
the covariance matrix $\Sigma$ is non degenerate, the dominant spectral value
$k(\lambda)$ satisfies: $r(P_\lambda) = |k(\lambda)| < 1$, for
$\lambda$ small and $\not = 0$. Hence $R_\mu \cap W = \{0\}$ for
some neighborhood $W$ of 0, i.e.  $R_\mu$ is a discrete subgroup of
$\RR^d$. \eop

\vskip 3mm {\bf Remark.} If $c_a(x)$ takes values in $\ZZ^d$, the
previous results have an analogue if we replace the space $X \times
\RR^d$ by $E = X \times \ZZ^d$. The character $\lambda \in V$ should
be replaced by a character of $\ZZ^d$, i.e. $\lambda \in \TT^d$, and
the Lebesgue measure $\ell$ by the counting measure on $\ZZ^d$. This
will be used in \ref{oncoverings} below.

\vskip 3mm The following corollary makes explicit the result in Theorem
\ref{recurr} for a functional $c(x)$ of a Markov chain.

\vskip 3mm \begin{cor} \label{cor2} Assume that $(X, \nu, \Gamma,
\mu)$ satisfies property (SG), that $c_a(x) = c(x)$ does not
depend on $a \in A$, and that the $\RR^d$-valued function $c$ is bounded and
satisfies $\int c(x) \ d\nu(x) = 0$. Moreover, assume that $A \subset \Gamma$
has a symmetric subset $B$ such that $B^2$ acts ergodically on $(X, \nu)$.
Then we have: \hfill \break 1) if $d \leq 2$, we have $\PP\times \nu$-a.e.
$\liminf_{n \to \infty} \|S_n(\omega,x) \| = 0$;

2) if the measure $c(\nu)$ is not supported on a coset of a proper
closed subgroup of $\RR^d$, \hfill \break - for $d \leq 2$, $\tilde
\sigma$ is ergodic with respect to $\mu \times \nu \times \ell$,
\hfill \break - for $d \geq 3$, the local limit theorem is valid for
$(S_n(\omega,x))$ and  $\lim_{n \to \infty} \|S_n(\omega, x)\| =
+\infty$, $\PP \times \nu$-a.e.
\end{cor} \proof \ The result follows from Lemma \ref{cnotdep}
and Theorem \ref{recurr}. \eop

\vskip 3mm The arguments in the proof of the proposition give also
the following corollary, which is a direct consequence of the main
result of \cite{JoSc87} (see also Theorem 6.3 in \cite{FuSh99}).

\vskip 3mm \begin{cor} \label{cor312} Assume $\supp (\mu)$ is
finite, generates $\Gamma$ and the representation $\rho_0$ of
$\Gamma$ in $L_0^2(X)$ does not contain weakly $\Id_\Gamma$. Let
$\Gamma^*$ be the group of characters of $\Gamma$ and $\Gamma_X^*$
be the subset of elements of $\Gamma^*$ contained in the natural
representation of $\Gamma$ in $L^2(X, \nu)$. Then $\Gamma_X^*$ is a
finite subgroup of $\Gamma^*$. The measure $\mu$ satisfies property
(SG) if and only if $\supp (\mu)$ is not contained in a coset of the
subgroup $\ker \chi$ for some $\chi \in \Gamma_X^*$, $\chi \not =
1$. In particular, if $(\supp (\mu))^k$ generates $\Gamma$ for any
$k
> 0$, then $(SG)$ is satisfied.
\end{cor} \proof \ Let $\rho$ be the natural representation of
$\Gamma$ in $L^2(X, \nu)$. For a given $\chi \in \Gamma_X^*$, there
exists $\varphi \in L^2(X, \nu)$ with $\varphi(\gamma x) =
\chi(\gamma) \varphi(x)$, $\|\varphi\|_2 = 1$. The ergodicity of the
action of $\Gamma$ on $X$ implies that $\varphi$ is uniquely defined
up to a scalar, with $|\varphi(x)| = 1$, $\nu$-a.e.

Also it is clear that $\Gamma_X^*$ is a subgroup of $\Gamma^*$. To
obtain that $\Gamma_X^*$ is closed in $\Gamma^*$, we note that if
$\chi \in \Gamma^*$ satisfies for some sequence $(\varphi_n)$ with
$|\varphi_n(x)| = 1$, $\chi(\gamma) = \lim_n \varphi_n(\gamma x) /
\varphi_n(x)$, then $\chi \in \Gamma_X^*$. As in the proof of the
proposition this follows from Proposition 2.3 of \cite{Sc98}, since
the subgroup of $\TT$-valued coboundaries of $(\Gamma, X, \nu)$ is
closed in the group of cocycles endowed with the topology of
convergence in measure.

In order to show that each element of $\Gamma_X^*$ has finite order,
we observe that, using \cite{JoSc87}, $(\Gamma, X, \nu)$ has no non
atomic $\ZZ$-factor up to orbit equivalence. Hence, for every $\chi
\in \Gamma_X^*$ and some $n \in \NN^*$ one has $\chi^n = 1$. Since
$\Gamma$ is finitely generated, $\Gamma_X^*$ is a closed subgroup of
a torus. Therefore $\Gamma_X^*$ is finite.

If $\supp (\mu)$ is contained in the coset $\{\gamma \in \Gamma:
\chi(\gamma) = c\}$ and $\varphi(\gamma x) = \chi(\gamma) \,
\varphi(x)$ with $\varphi \in L^2(X)$, one has:
$$P\varphi(x) = \sum_{a \in \supp (\mu)} \varphi(ax) \, \mu(a) = c
\, \varphi(x),$$ hence $\mu$ does not satisfy (SG).

Conversely, if $\mu$ does not satisfy (SG), then for some $c$ of
modulus 1 and a sequence $(\varphi_n)$ in $L^2(X)$ with
$\|\varphi_n\|= 1$, we have $\lim_n \|P \varphi_n - c \varphi_n\|_2=
0$. As in the proof of the proposition, we can use the condition
that $\rho_0$ does not contain weakly $\Id_\Gamma$ to get that $P -
I$ is invertible on $L_0^2(X)$ and obtain that $\lim_n \|
|\varphi_n| - 1\|_2 = 0$. Then, writing $\varphi_n(x) =
|\varphi_n(x)| \, e^{i \alpha_n(x)}$, we get that for a subsequence
$(n_k)$ of integers, $\lim_k e^{i(\alpha_{n_k}(x) -
\alpha_{n_k}(\gamma x))} = c$.

Then there exists $\chi \in \Gamma^*$ with $\chi(a) = c$ for every
$a \in A$, and for every $\gamma \in \Gamma$, $\chi(\gamma) = \lim_k
e^{i(\alpha_{n_k}(x) - \alpha_{n_k}(\gamma x))}$. From above $\chi
\in \Gamma_X^*$. The condition $\chi(a) = c$ for every $a \in A$
implies $\supp (\mu) \subset \{\gamma \in \Gamma: \chi(\gamma) =
c\}$. Hence the result.

For the last assertion, we observe that, if $\chi \in \Gamma_X^*
\setminus \{1\}$ satisfies $\chi(a) = c$ for some $c \in \TT$ and
every $a \in \supp (\mu)$, then, for some $ k \in \NN^*$, $\chi^k(a)
= c^k = 1$. Then any $\gamma \in (\supp (\mu))^k$ satisfies
$\chi(\gamma) = 1$. Since $(\supp (\mu))^k$ generates $\Gamma$, we
get $\chi = 1$, which contradicts the hypothesis on $\chi$. \eop

 \vskip 5mm
\section{Examples \label{examples}}

\subsection{Random walk in random scenery \label{randScenery}}
As an example corresponding to Corollary \ref{cor2}, we consider a
group $\Gamma$, a probability measure $\mu$ on $\Gamma$ such that
$A:=\supp(\mu)$ is symmetric and $(\supp(\mu))^2$ generates $\Gamma$
as a group. We denote by $(\Sigma_n(\omega), \omega \in A^\ZZ)$ the
left random walk on $\Gamma$ defined by $\mu$ and we consider the
visits of $\Sigma_n(\omega)$ to a random scenery on $\Gamma$.

Such a random scenery is defined by a finite set $C \subset \RR^d$,
a probability measure $\eta$ on $C$ with $\supp(\eta) = C$ and
$\sum_{v \in C} \eta(v) \, v = 0$. To each $\gamma \in \Gamma$,
one associates a random variable $x_\gamma$ with values in $C$.
The variables $x_\gamma$ are assumed to be i.i.d. with law $\eta$.

The scenery defines a point $x =(x_\gamma)_{\gamma \in \Gamma}$ of
the Bernoulli scheme $X = C^\Gamma$ endowed with the product measure
$\mu^{\otimes \Gamma}$ and $\Gamma$ acts on $C^\Gamma$ by left
translations: if $x = (x_\gamma, \gamma \in \Gamma)$, then $ax =
(x_{a\gamma}, \gamma \in \Gamma)$. If we define $f(x) = x_e \in C$,
the cumulated scenery is given by $S_n(\omega,x) = \sum_{k=1}^n
f(a_k(\omega)...a_1(\omega) x)$. One can give the following
interpretation: the random walker collects the quantity $x_\gamma$
when visiting the site $\gamma$ and his "cumulated gain" at time $n$
along the path defined by $\omega$ is $S_n(\omega, x)$.

The probability measure $\nu =\eta^{\otimes \Gamma}$ is $\Gamma$-invariant,
ergodic, and
$$\int f(x) d\eta^{\otimes \Gamma} (x) = \sum_{v \in C} \eta(v) v = 0.$$
The transformation $\sigma_1$ on $A^\ZZ \times C^\Gamma= \hat Y$
is given by $\sigma_1(\omega, x) = (\theta \omega, a_1(\omega) x)$.
Since $A$ is symmetric, $\sigma_1$ can be seen as a "$T -T^{-1}$" transformation
(cf. \cite{Kal82}).

We consider the Markov operator $P$ on $X$ associated with $\mu$ and
its restriction $\Pi_0$ to $L_0^2(X)$. It is well known (see \cite{BeHaVa08},
Ex E45, p. 394) that the action of $\Gamma$ on $L_0^2(X)$ decomposes as a direct sum of
tensor products of the regular representation in $\ell^2(\Gamma)$. A
typical summand is $\otimes_1^k \ell^2(\Gamma)$ and if $r_0(\mu)$ is
the spectral radius of the convolution operator by $\mu$ in
$\ell^2(\Gamma)$, we have $r(\Pi_0) = \sup_{k \geq 1}(r_0(\mu))^k =
r_0(\mu)$.

Assume that $\Gamma$ is non amenable. Then we have $r_0(\mu) < 1$ (see \cite{Ke59}),
hence property (SG) is satisfied. If we assume that $C \subset
\RR^d$ is not supported on a coset of a closed subgroup of $\RR^d$,
the hypothesis and the conclusions 2 and 3 of Corollary \ref{cor2}
are valid. Hence, with the above notations, it follows:

\begin{prop} Let $\Gamma$ be a non amenable group, $\mu$ a
probability measure on $\Gamma$ such that $\supp(\mu)$ is symmetric
and $(\supp(\mu))^2$ generates $\Gamma$, $\Sigma_n(\omega)$ the
corresponding random walk on $\Gamma$. We assume that $\Gamma$ is
endowed with an $\RR^d$-valued random scenery with law $\eta$, that
$C = \supp(\eta)$ is finite with $\sum_{v \in C} v \eta(v) = 0$, and
$\supp(\eta)$ is not contained in a coset of a closed subgroup of
$\RR^d$. We denote by $S_n(\omega, x)$ the accumulated scenery along
the random walk and by $\tilde \sigma$ the transformation on $\Omega
\times C^\Gamma \times \RR^d$ defined with $f(x) = x_e$ by
$$ \tilde \sigma(\omega, x, t) =(\sigma \omega, a_1(\omega) x, t + f(x)).$$
Then, the convergence of
${1 \over \sqrt n} S_n(\omega, x)$ to a non degenerate normal law is valid.
If $d \leq 2$, $\tilde \sigma$ is ergodic and $(S_n)$ is
recurrent with respect to $\mu \times \nu \times \ell$. If $d \geq
3$, $\mu^{\otimes \ZZ} \times \nu$-a.e.,
$$\lim_n \|S_n(\omega, x)\| = + \infty.$$
\end{prop}

\begin{rem} {
The above result should be compared with the case $\Gamma$ amenable.
For $\Gamma = \ZZ$, ${1 \over n^{3/2}} S_n(\omega, x)$ converges in law
towards a non Gaussian law  (\cite{KeSp79}, \cite{LeB06}).

Here, due to the strong transience properties of $\Gamma$, $S_n(\omega, x)$
behaves qualitatively like a sum of i.i.d. random variables.
Let us consider $\Gamma = \ZZ^m$, for $m$ large.

Using independence of the random variables $(x_\gamma, \gamma \in \Gamma)$,
we see that
$$\|\sum_{k= 0}^\infty P^k f\|_2^2 = \|f\|_2^2 \sum_{\gamma \in \Gamma} (\pi(\gamma))^2,$$
where $\pi = \sum_{k= 0}^\infty \mu^k$ is the potential of $\mu$ on
$\ZZ^m$. If $m \geq 5$, it is known that $\sum_{\gamma \in \Gamma}
(\pi(\gamma))^2 < \infty$ (see for example \cite{Uc07}), hence $
\varphi = \sum_{k= 0}^\infty P^k f$ is finite $\eta^{\otimes
\Gamma}$-a.e. and defines an element of $L_0^2(X)$ which satisfies
$(I - P) \varphi = f$. This implies the convergence in law of
${1\over \sqrt n} S_n(\omega, x)$ to a non degenerate normal law
(\cite{{GoLi78}}). }\end{rem}

\vskip 3mm
\subsection{Random walks on extensions of tori}

Now we present a special case where Condition $(\AP)$ can be
checked.

Here the $2d$-dimensional torus $\TT^d$ is identified with
$[-{1\over 2}, {1\over 2}[^{2d}$ and $\{x\}$ denotes the point of
$[-{1\over 2}, {1\over 2}[^{2d}$ corresponding to $x \in \TT^{2d}$.

\begin{prop} Let $\mu$ be a probability measure on $Sp(2d, \ZZ)$
acting by automorphisms on $\TT^{2d}$ and let $\Gamma$ be the
subgroup generated by $\supp(\mu)$. Assume that $\Gamma$ acts
$\QQ$-irreducibly on $\RR^{2d}$ and $\Gamma$ is not virtually
abelian. Let $\nu$ be the Lebesgue measure on $\TT^{2d}$.
Then, with the notations of Section \ref{appli}, we consider the
transformation $\tilde \sigma$ on $\Omega \times \TT^{2d} \times
\RR^{2d}$ defined by
$$ \tilde \sigma(\omega, x, v) =(\sigma \omega, a_1 x, v + \{x\}).$$
Let $S_n(\omega, x) = \sum_{k=1}^n \{a_k...a_1 x\}$. Then, if $d
\leq 2$, $\tilde \sigma$ is ergodic and $(S_n)$ is recurrent with
respect to $\mu^{\otimes \ZZ} \times \nu \times \ell$. If $d \geq 3$,
$\mu^{\otimes \ZZ} \times \nu$-a.e. $\lim_n \|S_n(\omega, x)\| = + \infty$.
\end{prop}

Since $x \to \{x\}$ is bounded and $\int \{x\} \, d\nu (x) =
0$, the proposition is a direct consequence of Proposition \ref{P0},
Theorem \ref{recurr} and the following lemma.

\begin{lem} Let $\mu$ be a probability measure on $Sp(2d, \ZZ)$ and let
$\Gamma$ be the subgroup generated by $\supp(\mu)$. For $\lambda \in
\RR^{2d}$ let $P_\lambda$ be the operator on $L^2(\TT^{2d})$ defined
by $$P_\lambda \varphi(x) = \sum_a e^{i\langle \lambda,
\{x\}\rangle} \varphi(ax) \, \mu(a).$$ Then, if $\Gamma$ acts
$\QQ$-irreducibly on $\RR^{2d}$ and $\Gamma$ is not virtually
abelian, we have $r(P_\lambda) < 1$, for $\lambda \in \RR^{2d}
\setminus \{0\}$. In particular (SG) and (AP) are valid.

\end{lem} \proof \ We will use as an auxiliary tool the Heisenberg group $H_{2d+1}$ and
its automorphism group $Sp(2d, \RR) \ltimes \RR^{2d}$. The group
$H_{2d+1}$ has a one dimensional center $C$ isomorphic to $\RR$ and
a lattice $\Delta$ such that $\Delta \cap C$ is isomorphic to $\ZZ$,
and $\Delta / \Delta \cap C$ is isomorphic to $\ZZ^{2d}$.

Let $\hat X$ be the corresponding manifold $H_{2d+1}/ \Delta$. Up to
a set of 0 measure, we can represent $\hat X$ as $\TT^{2d} \times
\TT^1$ and $\hat x \in \tilde X$ as $\hat x = (x,z)$, with $x \in
\TT^{2d}$, $z \in \TT^1 = \RR/\ZZ$.

The action of an element $g$ of $Sp(2d) \ltimes \RR^{2d}$ on
$H_{2d+1}$ can be represented as a matrix $g =\begin{pmatrix} (a) &
0 \\ u & 1 \end{pmatrix}$ acting on $\RR^{2d} \times \R$, where $a
\in Sp(2d, \RR)$, $u \in \RR^{2d}$. If $g$ preserves $\Delta$, we
have $a \in Sp(2d, \ZZ)$, $u \in \ZZ^{2d}$ and the action of $g$ on
$\hat X$ is given by $g(x,z) = (ax, z +[u,x])$, with $a \in Sp(2d,
\ZZ)$, $x = (x_1, x_2) \in \TT^{2d}$, $u = (u_1, u_2), u' = (-u_2,
u_1) \in \ZZ^{2d}$ and
$$[u, x] = \{\langle u_1, x_2 \rangle - \langle u_2, x_1
\rangle\} = \{\langle u', x \rangle\}  \in \TT^{2d}.$$

For a fixed $u \in \ZZ^{2d}$, we associate to $a \in Sp(2d, \ZZ)$
the element $g = \hat a$ of $Sp(2d, \ZZ) \ltimes \ZZ^{2d}$ with
components $a \in Sp(2d, \ZZ)$ and $u \in \ZZ^{2d}$.

We denote by $\hat \mu$ the push-forward of $\mu$ by the map $a \to
\hat a$. We denote by $\hat \Gamma$ the group generated by $\supp
(\hat \mu)$ and we consider the convolution action of $\hat \mu$ on
$L^2(\hat X)$. On functions $\psi$ of the form $\psi(x,z) =
\varphi(x) \, e^{2 i \pi z}$, with $\varphi \in L^2(\TT^{2d})$, the
action of $\hat \mu$ is given by
$$\hat P_u \varphi (x) = \sum_a \varphi(ax) \, e^{2 \pi i [u, x]} \,
\mu(a) = \sum_a \varphi(ax) \, e^{2 \pi i \langle u', x\rangle} \,
\mu(a) =  P_{2 \pi u'} \, \varphi (x).$$

We denote by $\|\mu^n\|_2$ (resp. $\|{\overline \mu}^n\|_2$) the
norm of the convolution operator by $\mu^n$ on $\ell^2(\Gamma)$
(resp. by ${\overline \mu}^n$ on $\ell^2(\ZZ^{2d} \setminus
\{0\})$). We write
\begin{eqnarray*}
r_0(\mu) &=& \lim_n \|\mu^n\|_2^{1/n}, \ r_0(\overline \mu) =\lim_n \|\overline \mu^n\|_2^{1/n}, \\
r_n &=& \sup(\|\mu^n\|_2, \|\overline \mu^n\|_2^{1/2d+2}).
\end{eqnarray*}

Now we observe that, for $\lambda, \lambda' \in \RR^{2d}$, we have
$$|(P_\lambda - P_{\lambda'}) \, \varphi| \leq {1\over 2} |\lambda - \lambda'|
\, P|\varphi|,$$ and $\|P_\lambda - P_{\lambda'}\| \leq {1\over 2}
|\lambda - \lambda'|$ since $P$ is a contraction.

On the other hand, if $\lambda' = 2\pi u'$, $u' \in \ZZ^d$, we have for any $n \in
\NN$, using \cite{BeHe11} (Theorem 3), $\|P_{\lambda'}^n \| \leq
r_n$.

We observe also that the qualitative result $r(P_{\lambda'}) < 1$
follows from the implication $b) \Rightarrow a)$ of Proposition
\ref{P0} and the remark following it, since $\Gamma$ acts
ergodically on $\TT^{2d}$, and therefore $\tilde \Gamma $ acts
ergodically on $\tilde X$ (cf. section 1).

The hypothesis on $\Gamma$ implies its non amenability (See
\cite{BeGu11}, Corollary 6), hence (see \cite{Ke59}) the spectral
radius $r_0(\mu)$ of the convolution operator on $\ell^2(\Gamma)$
defined by $\mu$ satisfies $r_0(\mu) < 1$. Also from \cite{BeGu11},
Corollary 6, we have $r_0(\overline \mu) = r(\Pi_0) < 1$.

We bound $\|P_{\lambda}^n \|$ as follows: we have $P_{\lambda'}^n -
P_{\lambda}^n = \sum_{k=0}^{n-1} P_{\lambda'}^k (P_{\lambda'} -
P_{\lambda}) P_{\lambda}^{n-k-1}$. Since $P_\lambda$ is a
contraction on $L^2(\TT^{2d})$, we have:
\begin{eqnarray*}
\|P_{\lambda'}^n - P_{\lambda}^n\| \leq \sum_{k=0}^{n-1}
\|P_{\lambda'}^k\|\, \|P_\lambda - P_{\lambda'}\| \leq \|P_\lambda -
P_{2\pi u'}\| \sum_{k=0}^{n-1} r_k.
\end{eqnarray*}
Hence, $\|P_\lambda^n\| \leq c \|\lambda - 2 \pi u'\| +
r_n$, with $c = {1\over 2} \sum_{k=1}^\infty r_k$,
which is finite since $r_0(\mu) < 1$, $r_0(\overline\mu) < 1$.

Since $\lim_n r_n = 0$, in order to show that $r(P_\lambda)
< 1$, i.e., $\|P_\lambda^n\| < 1$ for some $n > 0$, it suffices to
find $u' \in \ZZ^{2d}$ such that $c \|\lambda - 2 \pi u'\| < 1$.
This is possible at least for a multiple of $\lambda$: one can find
$k \in \NN$, $k \not = 0$, and $u' \in \ZZ^{2d}$ such that $\| k
\lambda - 2 \pi u' \| < c^{-1}$.

Now, if $r(P_\lambda) = 1$, one has also from Corollary
\ref{discSubG} that, for any $k \in \ZZ$, $r(P_{k\lambda}) = 1$.
From above this is impossible; hence $r(P_\lambda) < 1$. \eop

\goodbreak \subsection{Random walks on coverings\label{oncoverings}}

Let $G$ be a Lie group, $H$ a closed subgroup such that $G/H$ has a $G$-invariant measure $m$.
If $\mu$ is a probability measure on $G$, we consider the random walk on $G/H$
defined by $\mu$, and the corresponding skew product $\tilde \sigma$ on
$G^\ZZ \times G/H$ endowed with the measure $\mu^{\otimes \ZZ} \times m$.
Then one can ask for the ergodicity of such a skew product and its stochastic properties.
If $H$ is a normal subgroup  of another group $L \subset G$ such that $G/L$ is compact,
$G/H$ is fibred over $G/L$ and one can use harmonic analysis on $G/L$ and $H/L$.

A special case of Proposition \ref{wlkcov} below corresponds to the abelian coverings
of compact Riemann surfaces of genus $g \geq 2$. In this case, $H$ is a subgroup
$\Delta'$ of a cocompact lattice $\Delta$ in $\SL(2, \RR)$
and $G/\Delta'$ can be seen as the unit tangent bundle of the covering.

\begin{prop} \label{wlkcov} Let $G$ be a simple non compact real Lie group of real rank 1,
$\mu$ a symmetric probability measure with finite support $A \subset G$ such
that the closed subgroup $G_\mu$ generated by $A$ is non amenable. Let
$\Delta$ be a co-compact lattice in $G$, $\Delta'$ a normal subgroup
such that $\Delta/\Delta' = \ZZ^d$, $m$ the Haar measure on
$G/\Delta'$.

Let $\tilde \sigma$ be the extended shift on $\Omega \times
G/\Delta'$ defined by $\tilde \sigma(\omega,y) = (\sigma \omega,
a_1(\omega) y)$ and write $\Sigma_n(\omega) = a_n...a_1 \in G$.

If $d \leq 2$, $\tilde \sigma$ is ergodic with respect to
$\mu^{\otimes \ZZ} \times m$. If $d \geq 3$, we have $\mu^{\otimes
\ZZ} \times m$-a.e. $\lim_n \Sigma_n(\omega) y = +\infty$.
\end{prop} \proof \ \ Since $\Delta'$ is normal in $\Delta$, the
group $\Lambda = \Delta/\Delta' \sim \ZZ^d$ acts by right
translations on $G/\Delta'$ and this action of $\Lambda$ commutes
with the left action of $G$.

The $G$-space $G/\Delta'$ can be written as $X \times \Lambda$ where
$X \subset G/\Delta'$ is a Borel relatively compact fundamental
domain of $\Lambda$ in $G/\Delta'$. We will denote by $\overline y$
the projection of $y \in G/\Delta'$ on $X$ identified with
$G/\Delta$, by $\overline m$ the Haar measure on $G/\Delta$, and by
$(g, x) \to g. x$ the natural action of $g \in G$ on an element $x$
of the fundamental domain $X$

Let $z(y)$ be the $\Lambda$-valued Borel function on $G/\Delta'$
defined by $y = \overline y z(y)$. Then the $G$-action on $X \times
\Lambda$ can be written as $g(x,t) = (g.x, t + z(gx) )$ where the
group $\Lambda = \ZZ^d$ is written additively.

For $g \in G$ and $x \in X$, writing $Z(g,x) :=
z(g x)$, we obtain a cocycle:
$$Z(g_2 g_1,x) = Z(g_2,g_1.x) + Z(g_1,x).$$
Actually the cocycle relation is valid in restriction to $\Gamma$.

\vskip 3mm Since $G$ is simple and $G_\mu$ is non amenable, we know (\cite{FuSh99},
Theorem 6.11) that the convolution operator $\Pi_0$ on $X = G/\Delta$
defined by $\mu$ has a spectral radius $r(\Pi_0) < 1$ on $L_0^2(X)$.
On the other hand, if for any $a \in \supp(\mu)$, $x \in X$, we
write $c_a(x) = z(ax)= Z(a,x)$ and $\tilde a(x,t) = (a.x, t +
c_a(x))$, we are in the situation of Section \ref{appli}.

In order to verify this, we observe that, since $X$ is relatively
compact and $\supp(\mu)$ is finite, the functions $c_a(x)$ are
uniformly bounded. Furthermore, the cocycle relation for $Z(g,x)$
gives for any $g \in G$, $x \in X$: $Z(g^{-1},x) + Z(g,g^{-1}.x) =
0$; hence $\int (Z(g^{-1}, x) + Z(g,x)) \, d\overline m(x) = 0$.
Since $\mu$ is symmetric, we have the centering condition: $\int
c_a(x) \, d\overline m (x) \ d\mu(a) = 0$.

For any character $\lambda \in \Lambda^*$, any $\varphi \in L^2(X)$,
the formula $\rho_\lambda(g) \varphi(x) = e^{i\langle \lambda,
z(g^{-1}x)\rangle} \varphi(g^{-1} . x)$ defines a unitary one-dimensional
representation of $G$, hence of the group generated by $\supp(\mu)$, since
$Z(g,x)$ satisfies the cocycle relation.

Hence, using Proposition \ref{P0} and Theorem \ref{recurr}, the
proof will be finished if we show that $r(\rho_\lambda(\mu)) < 1$, for
$\lambda \not = 0$.

Since $G_\mu$ is non amenable and $G$ is simple, the result will
follow from Theorem C, part 2 of \cite{Sh00}, if we can show that
$\rho_\lambda$ does not contain weakly the representation $\Id_G$.
By definition, $\rho_\lambda$ is the induced representation to $G$
of the representation $\lambda_\Delta$ of $\Delta$ defined by the
character $\lambda$. Clearly, if $\lambda \not = 0$,
$\lambda_\Delta$ does not contain weakly $\Id_\Delta$. Since
$G/\Delta$ has a finite $G$-invariant measure, it follows from
Proposition 1.11b, p. 113 of \cite{Ma91} that $\rho_\lambda$ does
not contain weakly $\Id_G$. \eop

\vskip 2mm \goodbreak \subsection{Random walks on motion groups} \label{motGroup}

Let $G$ be the motion group $SU(d) \ltimes \CC^d$, $d \geq 2$. Write
$X =SU(d)$, $\nu$ for the Haar measure on $X$, $V =\CC^d$. We
identify a vector in $V$ with the corresponding translation in $G$
and we write $G = X \, V$. Let $\Gamma \subset SU(d)$ be a dense
subgroup with property (SG) and $A$ a finite generating set of
$\Gamma$. As mentioned is Section 2, such groups exists if $d \geq
2$. To each $a \in A$ we associate $\tilde a \in G$ with $\tilde a =
a \, \tau_a$, where $\tau_a \in V$. We consider a probability
measure $\mu$ on $A$ with $\supp(\mu) = A$ and we denote by $\tilde
\mu$ its push-forward on $\tilde A:= \{\tilde a, a \in A\}$.

In contrast to the above examples the main role here will be played by
$\tilde \Gamma$, the subgroup of $G$ generated by $\tilde A$.
Let us consider the convolutions $\tilde \mu^n, n \in \N$,
on $G$ and the natural affine action of $G$ on $V$.

\begin{proposition} \label{46} Assume that $\Gamma \subset SU(d)$ is such that
the natural representation of $\Gamma$ in $L_0^2(SU(d))$ does not
contain weakly $\Id_\Gamma$ and the affine action of $\tilde A$ on
$V$ has no fixed point. Then there exists $c
> 0$ such that for any continuous function $f$ with compact support
on $G$, $\lim_n \tilde \mu^n (f) \, n^{d} = c(\nu \otimes \ell)(f)$.
In particular, for any $f, f'$ continuous non negative functions on
$G$ with compact support, we have:
$$\lim_n {\tilde \mu^n(f) \over \tilde \mu^n(f')}
= { \int f(g) \, dg  \over \int f'(g) \, dg }.$$

Furthermore the convolution equation $\tilde \mu * f = f$ on $G$,
with $f \in L^\infty(\nu \otimes \ell)$, has only constant
solutions. \end{proposition} \proof \ We will use the results of
Section 3; the link with Section 3 is as follows. The maps $\tilde
a$ on $X \times V$ are defined here as left multiplication on $G = X
V$ by $a \tau_a$:
$$\tilde a (x, v) = \tilde a (xv) = ( ax, v + x^{-1} (\tau_a)),$$
where $x^{-1} (\tau_a)$ is the vector obtained from $\tau_a$ by the
linear action of $x$.

Hence the action of $A$ on $X$ is by left multiplication on the
group $SU(d)$ and $c_a(x) = x^{-1}(\tau_a)$. The centering condition
$\int c_a(x) \, d\mu(a) \, d\nu(x) = 0$ is valid here, since it
reduces to $\int x^{-1}(\tau_a) \, d\nu(x) = 0$, which is a
consequence of the fact that this integral is the barycenter of the
sphere $SU(d)\tau_a$ of center 0 and radius $\|\tau_a\|$, hence is
equal to 0.

Then the action of $\tilde \Gamma \subset G$ on $X \times V$ is by left multiplication
on $G = XV$. This action is part of the action of $G$ on itself by left translation.

Let us fix some notations. For $x \in X$, $v \in V$, with the above notations,
$x(v)$ corresponds to the element $x v x^{-1}$ of $G$.
We observe that, if $g = x_g \tau_g$ and $h = x_h \tau_h$,
then $x_{gh} = x_g x_h$, $\tau_{gh} = x_h^{-1}(\tau_g) + \tau_h$.

Therefore $(g,x) \to x^{-1}(\tau_g)$ is a $V$-valued cocycle on $G \times X$, where the action
of $G$ on $X$ is given by $(g,x) \to x_g x$:
$$x^{-1}(\tau_{gh}) = (x_hx)^{-1}(\tau_g) + x^{-1}(\tau_h).$$
It follows that for $\tilde \gamma \in \tilde \Gamma$, $c(\tilde \gamma,x)$ as defined in Section 3 is equal to
$x^{-1}(\tau_{\tilde \gamma})$ and is the restriction to $\tilde \Gamma \times X$
of the cocycle on $G \times X$ given by $c(g,x) = x^{-1}(\tau_g)$.

We show now that the closure $H$ of $\tilde \Gamma$ is equal to $G$.
We observe that $H \cap V$ is a normal subgroup of $H$ and the
action by conjugation of $G$ on $V$ reduces to the linear action of
$G$.

Since $\Gamma$ is dense in $SU(d)$ and $W= H \cap V$ is
$\Gamma$-invariant, $W$ is a closed $SU(d)$-invariant subgroup of
$V$. Hence $W=\{0\}$ or $V$.

Suppose we are in the first case. Then the projection $H \to SU(d)$
is injective. In particular $\tilde \Gamma$ is isomorphic to
$\Gamma$. In connection with Section 3, we may observe that
$c(\tilde \gamma, x) = x^{-1}(\tau_{\tilde \gamma}) $ defines also a
cocycle on $\Gamma \times X$ since $\tau_{\tilde \gamma}$ depends
only on $\gamma$; hence $c(\tilde \gamma, x) = c(\gamma, x)$.

We will use the following lemma.
\begin{lem} Assume $H$ is a closed subgroup of $G = SU(d) \ltimes \CC^d$, $d \geq 2$, such that $H
\cap \CC^d = \{0\}$ and the projection of $H$ on $SU(d)$ is dense.
Then $H$ is conjugate to $SU(d)$.
\end{lem} \proof \  Let $\pi$ be the projection of $G$ onto $SU(d)$.
Observe that $\pi(H)$ is a Lie subgroup of $SU(d)$ isomorphic to
$H$. Also $\pi(H)$ contains a finitely generated countable subgroup
$\Delta$ which is dense in $\pi(H)$, hence in $SU(d)$. Then $\Delta$
is non amenable since otherwise, using \cite{Ti72}, $\Delta$ would
have a polycyclic subgroup $\Delta_0$ with finite index. Then the
closure of $\Delta_0$ would be solvable and equal to $SU(d)$, which
is impossible since $d \geq 2$.

Let $H_0$ be the connected component of identity in $H$ and observe
that $\pi(H_0)$ is normal in $\pi(H)$. It follows that the Lie
algebra of $\pi(H_0)$ is invariant under the adjoint action of
$\pi(H)$, hence invariant under the action of its closure $SU(d)$.
Then, using the exponential map, we see that $\pi(H_0)$ is a normal
Lie subgroup of $SU(d)$.

Since $SU(d)$ is a simple Lie group, we get $\pi(H_0) = \{e\}$ or
$\pi(H_0) = SU(d)$. In the first case, $H$ would be a discrete
subgroup of $G$, hence amenable like $G$. This imply that
$\pi^{-1}(\Delta) \subset H$ would be amenable. Hence $\Delta$
itself would be amenable which is a contradiction. Hence $\pi(H) =
SU(d)$ and $\pi$ is an isomorphism of $H$ onto $SU(d)$. In
particular $H$ is compact and its affine action on $V$ has a fixed
point $\tau \in V$. Hence $\tau^{-1} H \tau = SU(d)$. \eop

The existence of a fixed point for the affine action of $H$ on $V$,
as shown in the lemma, contradicts the hypothesis on $\tilde A$,
hence $W = V$. Since the projection of $H$ on $SU(d)$ is dense, we
get $H = G$.

Now we are going to apply Theorem \ref{recurr}, part 1b). For this
we have to verify $(\AP)$. If (AP) is not valid, there exists
$(\lambda, \theta) \in V \times \RR$ and $d(x)$ with $|d(x)| = 1$,
such that for any $a \in A$:
$$e^{i\langle \lambda, c(a,x) \rangle} = e^{i\theta} d(x a) /d(x).$$

As observed above, $c(a,x)$ extends to $G$ as the cocycle $c(g, x)$
which is equal to $x^{-1}(\tau_{g})$ on $(g, x) = (x_g  \tau_{g},
x)$. Then we have
$$e^{i\theta} = e^{i\langle \lambda, c(a,x) \rangle} d(x) / d(a.x),$$
and the right hand side is the restriction to $\tilde A\times X$
of the cocycle
$$c_\lambda(g,x) = e^{i \langle\lambda, c(g,x)\rangle}
d(x) /d(x_g x)$$ on $G\times X$.

This cocycle takes values $e^{i\theta}$ on $\tilde A$, hence its values are
also independent of $x$ on the group $\tilde \Gamma$.
Since $\tilde \Gamma$ is dense in $G$ and $c_\lambda$ is measurable on $G \times X$,
using the $L^2$ continuity of the translation,
it is also independent of $x$ on $G$, hence it defines a character on $G$.

Since $G$ has no non trivial character we get $e^{i\theta} = 1$.
Then we have, for any $\tilde \gamma \in \tilde \Gamma$ with $\tilde
\gamma = \gamma \tau_{\bar \gamma}$ and a.e. $x \in X$,
$$e^{i\langle \lambda, x^{-1}(\tau_{\tilde \gamma}) \rangle} = d(\gamma x) / d(x).$$
This means that the function on $G$ defined by $\psi(xv) =
e^{-i\langle \lambda, v\rangle}\, d(x)$ is invariant by left
translation by any element $\tilde \gamma \in \tilde \Gamma$. Since
$\tilde \Gamma$ is dense in $G$, hence ergodic on $G$, $\psi$ is
constant, i.e. $\lambda = 0$, $d = 1$. It follows that (AP) is
valid. Hence the result.

Since $(\AP)$ is valid, the last assertion is a consequence of 3) in
Theorem \ref{recurr}. \eop

\vskip 3mm There exists various possibilities for the geometry of
the subgroup $\tilde \Gamma$ inside $G$, as the following
proposition shows.

\begin{prop} With the above notations, assume that the finite set $A \subset
SU(d)$ generates a dense subgroup $\Gamma$ and the affine action of
$\hat A$ on $V$ has no fixed point. Then

1) If $\Gamma$ has property (T), then $\tilde \Gamma \cap V$ is
dense in $V$.

2) If $\Gamma$ is a free group, then $\tilde \Gamma \cap V = \{0\}$,
hence $\tilde \Gamma$ is a dense subgroup of $G$ isomorphic to
$\Gamma$.
\end{prop} \proof \ 1) We show using arguments as in the proof of
Proposition \ref{46} that $\tilde \Gamma$ is dense in $G$. We
observe that $\tilde \Gamma \cap V$ is a normal subgroup of $\tilde
\Gamma$ and the action by conjugation of $G$ on $V$ reduces to the
linear action of $SU(d)$.

Since $\Gamma$ is dense in $SU(d)$ and $\tilde \Gamma \cap V$ is
$\Gamma$-invariant, its closure $W$ is a closed $SU(d)$-invariant
subgroup of $V$. Hence $W=\{0\}$ or $V$.

Suppose $W=\{0\}$. Then the projection $\tilde \Gamma \to \Gamma $
is injective, hence $\tilde \Gamma$ is isomorphic to $\Gamma$ and
has property (T). We have also $c(\tilde \gamma, x) =
x^{-1}(\tau_{\tilde \gamma}) = x^{-1}(\tau_{\gamma}) = c(\gamma,
x)$. Then the cocycle $c(\tilde \gamma, x)$ from $\tilde \Gamma
\times X$ to the vector group $V$ is trivial (See Zimmer, p. 162),
hence $c(\tilde \gamma, x) = \varphi(\gamma x) - \varphi(x)$ for
some $\varphi \in L^2(X)$.

Also, from above, $\tau_{\gamma} = x(c(\tilde \gamma,x))$ does not
depend on $x$. Then, for every $\gamma \in \Gamma$,
$$\tau_{\gamma} = x(\varphi(\gamma x) - \varphi(x))
= \int x(\varphi(\gamma x) - \varphi(x)) \, d\nu(x) = \gamma^{-1}(w)
- w,$$ with $w := \int x(\varphi(x))\, d\nu(x) =
 \int \gamma x(\varphi(\gamma x)) \, d\nu(x)$.

It follows, for the affine action of $\tilde \gamma$ on $V$:
$$\tilde \gamma w = \gamma(w + \tau_\gamma) = \gamma(\gamma^{-1}(w))  = w.$$

This contradicts the hypothesis on $\tilde A$, hence $\tilde \Gamma
\cap V = \{0\}$ is not valid. Therefore $\tilde \Gamma \cap V$ is
dense in $V$. The fact that $\tilde \Gamma$ is dense in $G$ follows,
but was already proved in Proposition \ref{46}.

\vskip 3mm 2) We denote by $\pi$ the natural projection of $G$ onto
$SU(d)$, and we observe that $\pi(\tilde a) = a$ for any $a \in A$;
hence $\pi(\tilde \Gamma) = \Gamma$. Since $\Gamma$ is free it
follows that the restriction of $\pi$ to $\tilde \Gamma$ is an
isomorphism of $\tilde \Gamma$ onto $\Gamma$.

In particular, since $\pi(V)= \{0\}$ we have $\tilde \Gamma \cap V =
\{0\}$ and $\tilde \Gamma$ is free. The density of $\tilde \Gamma$
in $G$ has been shown in the proof of Proposition \ref{46}.

 \eop

\section{Questions}

1) In the situation of random walks in random scenery (example
\ref{randScenery}), with $\Gamma = \ZZ^m$, $m \geq 3$, is the local
limit theorem for $S_n(\omega, x) \in V$ valid ?

2) In the situation of motion groups (example \ref{motGroup}), for
$d \geq 2$, if $\tilde \Gamma \subset SU(d) \ltimes \CC^d$ and
$\Gamma \subset SU(d)$ is dense, is the local limit theorem for
$S_n(\omega, x) \in V$ still valid ?

What can be said about the equidistribution of the orbits of $\tilde \Gamma$ on $V$ ?

What are the bounded solutions of the equation $\tilde \mu
* f = f$, $f \in L^\infty(\nu \otimes \ell)$, on $G$.

 3) In the
above considerations the maps $a \in A$ are chosen with probability
$\mu(a)$ which does not depend on $x$ and the product space is
endowed with the product measure $\PP = \mu^{\otimes \NN^*}$. One
can extends this framework by choosing the maps $a \in A$ according
to a weight $\mu(x, a)$ depending on $x \in X$ and consider the
corresponding Markovian model. One can also replace the shift
invariant measure $\PP$ by a Gibbs measure.

A question is then the validity of  the results obtained above in
these more general situations.

\bibliographystyle{amsalpha}

\end{document}